\documentclass[12pt,a4paper,reqno]{amsart}
\usepackage{fullpage}
\usepackage{amsmath}
\usepackage{amsfonts}
\usepackage{amssymb}
\usepackage{amsthm}
\usepackage{mathrsfs}

\newtheorem{letterthm}{Theorem}

\newcommand{\ep}{\varepsilon}

\newcommand{\con}{\equiv}

\newcommand{\maps}{\rightarrow}

\newcommand{\al}{\alpha}
\newcommand{\be}{\beta}

\newcommand{\ga}{\gamma}

\newcommand{\sig}{\sigma}
\newcommand{\lam}{\lambda}

\newcommand{\beq}{\begin{equation}}
\newcommand{\eeq}{\end{equation}}

\begin{document}

\theoremstyle{plain}
\newtheorem{thm}{Theorem}[section]
\newtheorem{lemma}[thm]{Lemma}
\newtheorem{prop}[thm]{Proposition}
\newtheorem{question}[thm]{Question}
\theoremstyle{remark}
\newtheorem{remark}[thm]{Remark}

\def\beq{\begin{equation}}
\def\endeq{\end{equation}}

\newcommand{\N}{\mathbb{N}}
\newcommand{\R}{\mathbb{R}}
\newcommand{\Z}{\mathbb{Z}}
\newcommand{\C}{\mathbb{C}}
\newcommand{\Q}{\mathbb{Q}}

\newcommand{\Int}{\int\limits}
\newcommand{\DInt}{\Int_0^\infty}
\newcommand{\RInt}{\Int_{\R}}

\newcommand{\re}{\,\textrm{Re}\,}

\numberwithin{equation}{section}

\title[Polynomial Carleson operators along monomial curves]{Polynomial Carleson operators along monomial curves in the plane}
\date{\today}
\author[S. Guo, L. B. Pierce, J. Roos, and P.-L. Yung]{Shaoming Guo, Lillian B. Pierce, Joris Roos and Po-Lam Yung}
\address{}

\begin{abstract}
We prove $L^p$ bounds for   partial polynomial Carleson operators along  monomial curves $(t,t^m)$ in the plane $\R^2$ with a phase polynomial consisting of a single monomial. These operators are ``partial'' in the sense that we consider linearizing stopping-time functions that depend  on only one of the two ambient variables. A motivation for studying these partial operators is the curious feature that, despite their apparent limitations, for certain combinations of  curve and phase,  $L^2$ bounds for partial operators along curves imply the full strength of the $L^2$ bound for a one-dimensional Carleson operator, and for a quadratic Carleson operator.
Our methods, which can at present only treat certain combinations of curves and phases, in some cases adapt a $TT^*$ method to treat phases involving fractional monomials, and in other cases use a known vector-valued variant of the Carleson-Hunt theorem.
\end{abstract}

\maketitle

\section{Introduction}

\subsection{Historical background.} In 1966,  Carleson \cite{Carleson66} proved an $L^2$ bound for the Carleson operator
\beq
 \label{eq:carleson} f(x)\longmapsto \sup_{N\in\R} \Big| p.v.\int_\R f(x-t) e^{iN t} \frac{dt}t\Big|.
 \eeq
This provided the key step in proving almost everywhere convergence of Fourier series of $L^2$ functions and thereby resolved a conjecture of Luzin. The $L^p$  boundedness of the Carleson operator for $1<p<\infty$ was then shown by Hunt \cite{Hunt67}, and further proofs of Carleson's theorem were later given by Fefferman \cite{FeffermanCarleson} and Lacey and Thiele \cite{LaceyThieleCarleson}.

E. M. Stein suggested the following generalization: fix a natural number $d$ and consider the operator given by
\beq
\label{eqn:polynomialcarleson}
f(x)\longmapsto\sup_{P}\left|\int_{\R^n} f(x-y) e^{iP(y)} K(y) dy\right|,
\eeq
where $K$ is an appropriately chosen Calder\'on-Zygmund kernel and the supremum runs over all real-valued polynomials $P$ of degree at most $d$ in $n$ variables. Stein asked whether this \emph{polynomial Carleson operator} is bounded from $L^p$ to $L^p$ for $1<p<\infty$. 
Stein and Wainger \cite{SteinWainger01} used a $TT^*$ argument and certain oscillatory integral estimates of van der Corput type to obtain $L^p$ bounds for a variant of the operator \eqref{eqn:polynomialcarleson}, where the polynomial $P$ is restricted to the set of polynomials of degree at most $d$ that vanish to at least second order at the origin (so in particular, have no linear term; of course constant terms may be disregarded).
In dimension $n=1$, a positive answer to Stein's full question was provided by Lie \cite{Lie09}, \cite{Lie11}, who developed a sophisticated time-frequency approach.
In higher dimensions $n>1$, boundedness of the full polynomial Carleson operator  remains an open problem.

Pierce and Yung \cite{PierceYung15} have introduced a new aspect to the study of polynomial Carleson operators, by considering an operator that also features Radon-type behavior in the sense of integration along an appropriate hypersurface. More precisely, they considered the operator
\begin{align}
\label{eqn:carlesonpierceyung}
f(x,y) \longmapsto \sup_{P}
\left|\int_{\R^n} f(x-t,y-|t|^2) e^{iP(t)} K(t) dt\right|,
\end{align}
acting on functions $f$ on $\R^n\times\R$ where $n\ge 2$, $K$ is a Calder{\'o}n-Zygmund kernel, and the supremum runs over a suitable vector subspace of the space of all real-valued polynomials $P$ of degree at most $d$ in $n$ variables. In particular, this allowable subspace requires that the polynomials considered should omit linear as well as certain types of quadratic terms. The key result of \cite{PierceYung15} then proves $L^p$, $1<p<\infty$, bounds for this operator, via a method of proof based on square functions, $TT^*$ techniques in the spirit of Stein and Wainger \cite{SteinWainger01}, and certain refined van der Corput estimates. 
Notably, the method of \cite{PierceYung15} does not work in the planar case $n=1$, which is the main subject of the present paper. Our goal here is to prove bounds for a new class of polynomial Carleson operators along curves in the plane, and to demonstrate the curious feature that even \emph{partial} results for these new operators along curves are in some sense as strong as Carleson's original theorem (and its variants) in the purely one-dimensional setting.

\subsection{Statements of main results.}

Let $m,d$ be positive integers and $f$ a Schwartz function on $\R^2$. For $N\in\R$ let
\[ H_N f(x,y)=H^{m,d}_N f(x,y)=p.v. \RInt f(x-t, y-t^m) e^{iNt^d} \frac{dt}t.\]
The natural goal, in the spirit of Carleson operators, is to prove that for all $1<p<\infty$,
\beq\label{H_main}
 \left \| \sup_{N \in \R} |H_N f|  \right \|_{L^p(dxdy)} \leq C \|f\|_{L^p(dxdy)}. 
 \eeq
This would be analogous to the results of Stein and Wainger \cite{SteinWainger01} in the Radon-type context of (\ref{eqn:polynomialcarleson}).
We recall the useful strategy of linearization via a linearizing stopping-time function: we  define for an arbitrary measurable function $N(x,y): \R^2 \mapsto \R$ the operator $f \mapsto H_{N(x,y)}f(x,y)$. Then proving 
\[\| H_{N(x,y)}f(x,y)\|_{L^p(dxdy)} \leq C \|f\|_{L^p(dxdy)}\]
with a constant $C$ independent of the choice of the function $N$ is equivalent to proving (\ref{H_main}).

To prove this inequality appears to be out of reach of our current methods.
Recalling instead that a special case of \cite{SteinWainger01} already shows that  for any integer $d>1$ the operator
\beq\label{fop}
f(x) \longmapsto p.v. \RInt f(x-t) e^{iN(x) t^d} \frac{dt}t
\eeq
is bounded on $L^p(\mathbb{R})$ for $1 < p < \infty$, we are motivated to consider the case when we twist the operator (\ref{fop}) with an additional Radon transform, while preserving the dependence of the linearizing function $N$ on one variable only.

Thus for an arbitrary measurable function $N:\R\to\R$, we define our main operators of interest:
\beq
A_N^{m,d}f(x,y):=H^{m,d}_{N(x)} f(x,y)
\endeq 
and 
\beq
B_N^{m,d} f(x,y):=H^{m,d}_{N(y)} f(x,y).
\endeq 
Before turning to our main results, we briefly note that certain special cases of these operators may be treated immediately: namely, for $d \geq 1$ the operators $A^{1,d}_N$ and $B^{1,d}_N$ are bounded on $L^p(\mathbb{R}^2)$ for $1 < p < \infty$. Indeed even the operator $\sup_{N \in \mathbb{R}} |H^{1,d}_N f(x,y)|$ is bounded on $L^p(\mathbb{R}^2)$ for $1 < p < \infty$. This follows immediately by integrating Carleson's theorem for (\ref{eq:carleson}) (in the case $d = 1$), or the result of Stein and Wainger \cite{SteinWainger01} for (\ref{eqn:polynomialcarleson}) (in the case $d > 1$), along the straight lines of slope 1 in $\mathbb{R}^2,$ using Fubini's thoerem.

The remaining cases, with $m > 1$, are highly nontrivial. We formulate our main results as two theorems, which despite superficial similarities have quite different flavors, due to the differing symmetry groups of the involved operators (see Section \ref{sec_symmetry}).
Our first main result can be stated as follows.
\begin{thm}\label{2709theorem1.2}
\label{thm:main1}
Let \(N:\R\to\R\) be a measurable function and \(d,m> 1, d\not=m\) integers. Then for $1<p<\infty$,
\begin{align}
\label{eqn:mainNx}
 \Big\|A_N^{m,d} f\Big\|_{p} &\le C \|f\|_p, \\
\label{eqn:mainNy}
\Big\|B_N^{m,d} f\Big\|_{p} &\le C \|f\|_p,
\end{align}
 with the constant $0<C<\infty$ depending only on \(d,m,p\) and not on \(N,f\).
\end{thm}

Note that uniformity of \eqref{eqn:mainNx} in $N$ is tantamount to the estimate
\begin{align}\label{eqn:partialLp1}
\left\|\sup_{N\in\R} \| H_N f(x,y)\|_{L^p(dy)} \right\|_{L^p(dx)}\le C \|f\|_p.
\end{align}
Similarly, \eqref{eqn:mainNy} corresponds to
\begin{align}\label{eqn:partialLp2}
\left\|\sup_{N\in\R} \| H_N f(x,y)\|_{L^p(dx)} \right\|_{L^p(dy)}\le C \|f\|_p. 
\end{align}

Our proof of Theorem~\ref{2709theorem1.2} proceeds via van der Corput estimates, and does not depend on Carleson's theorem; this is in contrast to our second result, which we state as follows.

\begin{thm}\label{thm:main2}
Let \(N:\R\to\R\) be a measurable function. Then for $1<p<\infty$,
\begin{align}
\label{eqn:mainNx2}
 \Big\|A_N^{m,1} f\Big\|_{p} &\le C \|f\|_p, \quad \text{for any integer $m \geq 3$,}\\
\label{eqn:mainNy2}
\Big\|B_N^{m,m} f\Big\|_{p} &\le C \|f\|_p, \quad \text{for any integer $m \geq 2,$}
\end{align}
with the constant $C$ depending only on $m,p$ and not on $N,f$.
\end{thm}
A novel feature of our proof of Theorem \ref{thm:main2} is that we combine Carleson's theorem with $TT^*$ estimates in the spirit of Stein and Wainger. One surprising feature of our proof, compared to the original work \cite{SteinWainger01} is that these $TT^*$ estimates can handle certain cases of phase polynomials with a linear term (c.f. estimates \eqref{eqn:symmphasefct}--\eqref{2709ee3.26}).

\begin{remark}\label{remark_difficult}
One is led to ask what happens to the remaining nontrivial ($m >1$) cases that are not covered by Theorems~\ref{thm:main1} and \ref{thm:main2}, namely $A_N^{2,1}$, $A^{m,m}_N$ and $B^{m,1}_N$ where $m > 1$ is an integer. The key again lies in the symmetries of these operators: they are different from the symmetries of the operators in Theorems~\ref{thm:main1} and \ref{thm:main2}, and this points to why our current proofs do not apply in these situations.
Despite these difficulties, at least the $L^2$ bounds for all these problematic cases still follow from known Carleson theorems via partial Fourier transform and Plancherel's theorem; see Section \ref{sec_deductions2}. The full $L^p$ bounds remain an open problem in these cases. 

\end{remark}

\subsection{Consequences of bounding partial Carleson operators}\label{sec_consequences}
We now turn to the surprising feature that $L^2$ bounds for partial operators along curves  imply $L^2$ bounds for Carleson-type operators acting on functions of one variable. Here we summarize several deductions of this kind; proofs are given in Section \ref{sec_deductions}.

First, $L^2$ bounds for certain operators $A^{m,1}$ and $B^{m,m}$ are in some sense equivalent to an $L^2$ bound for Carleson's operator. More precisely, for any integer $m \geq 1$, the $L^2$ boundedness of $A^{m,1}_N$ implies the $L^2$ boundedness for the one-dimensional Carleson operator (\ref{eq:carleson}), by a Plancherel argument (see Section 6.1). In the other direction, we use the boundedness of the maximal truncated Carleson operator \eqref{define_Cstar} (itself dominated by the Carleson operator according to the inequality \eqref{CCM}) to prove the $L^2$ boundedness of $A^{m,1}$ for $m \geq 3$ in Theorem \ref{thm:main2}, while the $L^2$ bound for $A^{2,1}$ may be deduced from Carleson's theorem (see Section \ref{sec_deductions2}).

Similarly, for any \emph{odd} integer $m \geq 1$,  the $L^2$ boundedness of $B^{m,m}_N$ implies an $L^2$ bound for the one-dimensional Carleson operator (see Section \ref{sec_deductions1}), while in the other direction we use the maximal truncated Carleson operator to prove Theorem \ref{thm:main2}.

Of course, the most natural challenge in the setting of Carleson operators along curves in the plane is the quadratic Carleson operator along the parabola defined by
\begin{align}
\label{eqn:paraboliccarleson}
\mathscr{C}^{\text{par}} f(x,y)=\sup_{N\in\R^2}
\left|H_{N}^{\text{par}} f(x,y)\right|,
\end{align}
where for $f$ a Schwartz function on $\R^2$,
\[ H_{N}^{\text{par}} f(x,y)=p.v.\int_{\R} f(x-t,y-t^2) e^{iN_1t+iN_2t^2} \frac{dt}t. \]
This operator combines all the features that have proved troublesome in the study of \eqref{eqn:carlesonpierceyung} in \cite{PierceYung15}: apart from acting on functions in the  plane, the phase consists entirely of the problematic linear and quadratic terms. Assuming that $N_1,N_2:\R \to \R$ are  arbitrary measurable functions depending only on $x$, observe that for $N_1=0$ this gives our operator $A^{2,2}$ (which our present arguments cannot treat) and for $N_2=0$ it gives our problematic operator $A^{2,1}$ (which again our present arguments cannot treat). So we are quite far from knowing how to bound \eqref{eqn:paraboliccarleson}.

But in the spirit of studying partial versions of Carleson operators, we point out that even a partial estimate for $H_N^{\text{par}}$ of the form 
\begin{align}\label{partial_HPN}
\left\|\sup_{N\in\R^2} \|H_N^{\text{par}} f\|_{L^2(dy)} \right\|_{L^2(dx)}\le C\|f\|_2,
\end{align}
would suffice to imply an analogue over $\R$ of Lie's $L^2$ result on the quadratic Carleson operator \cite{Lie09}; see Section \ref{sec_deductions1} for details.
 These considerations indicate the interest in pursuing the partial Carleson operators we consider.

\section{Overview of the methods}

\subsection{Symmetries of our operators}\label{sec_symmetry}

To make precise the differences between Theorems \ref{thm:main1} and \ref{thm:main2}, we now characterize symmetries of the operators $A_N^{m,d}$ and $B_N^{m,d}$ as $m$ and $d$ vary.
First there is an anisotropic dilation symmetry. If we denote 
\beq
D_\lambda f(x,y)=f(\lambda x,\lambda^m y)
\endeq 
for $\lambda>0$, then
\[ D^{-1}_{\lambda} H^{m,d}_N D_\lambda = H^{m,d}_{\lambda^{-d} N}. \]
Second, due to the convolution structure, $H_N^{m,d}$ commutes with 
translations of the plane, for any $m,d$.

Third, the operators in Theorem \ref{thm:main2} additionally have certain modulation symmetries. Let 
\beq
M_{\xi,\zeta} f(x,y)=e^{ix\xi+iy\zeta} f(x,y)
\endeq 
for $\xi,\zeta\in\R$.  Then if $d=1$, we have
\beq\label{2809ee1.14} M_{\xi,0}^{-1} A^{m,1}_N M_{\xi,0} = A^{m,1}_{N-\xi}
\eeq for all $\xi \in \mathbb{R}$.
Similarly if $d=m$, we have
\beq\label{2809ee1.15} M_{0,\zeta}^{-1}  B^{m,m}_N M_{0,\zeta} = B^{m,m}_{N-\zeta}\eeq for all $\zeta \in \mathbb{R}$.
Simultaneous translation and modulation invariance is a characteristic property of the Carleson operator.
Hence we are led to use Carleson's theorem in parts of the proof of Theorem~\ref{thm:main2}.

Finally, we remark briefly that for $A^{2, 1}_N$, the modulation symmetries are more involved. The problem is that in addition to the modulation symmetry \eqref{2809ee1.14}, it also has a certain quadratic modulation symmetry. Let 
\beq
Q_{b}f(x, y)=e^{ibx^2} f(x, y).
\endeq  
Then 
\beq\label{2809ee1.32}
Q^{-1}_b M^{-1}_{0, b}A^{2,1}_N M_{0, b} Q_b=A^{2, 1}_{N-2bx}.
\endeq
Recall that for the operator $A^{2, 1}_N$, the linearizing function $N$ depends on the variable $x$. Thus, by $N-2bx$ we mean the function $x\mapsto N(x)-2bx$, also only depending on $x$. Moreover, notice that in \eqref{2809ee1.32}, the linear modulation acts on the $y$ variable, while the quadratic modulation acts on the $x$ variable. Hence there is a certain ``twist'' in this modulation symmetry. 

\subsection{Method of proof: Theorem \ref{thm:main1}}

We now sketch the proof of Theorem~\ref{thm:main1}. The strategy follows broadly that of Stein and Wainger, but the means of obtaining the key estimates is necessarily different. More precisely, we proceed by splitting the integral defining $A_N^{m,d}$ or $B_N^{m,d}$ into two parts, according to the size of the phase $Nt^d$: for $Nt^d$ sufficiently small, we compare the resulting operator to a maximal truncated Hilbert transform along a curve, and for $Nt^d$ large, we use $TT^*$ and van der Corput estimates to handle the operator that arises. It is in the treatment of this latter operator where we must assume the stopping time depends on one variable only, so that we may perform a Fourier transform in the free variable, along which the linearizing function is constant. This idea goes back to Coifman and El Kohen, who used it in the context of Hilbert transforms along vector fields (see the discussion in Bateman and Thiele \cite{BatemanThiele}).

Another important ingredient is a certain refinement of Theorem 1 of \cite{SteinWainger01}. The main novelty is that our core estimate, which we now record, allows us to consider phases with monomials of fractional exponents.

\begin{lemma}\label{lemma:main}
Fix real numbers $\alpha,\beta>0$, $\alpha\not=\beta$, $\alpha,\beta\not=1$. Let $\psi$ be smooth and supported on $[1,2]$. For $\lambda=(\lambda_1, \lambda_2)\in\R^2$ and $t > 0$, let 
\beq\label{Phi_dfn1}
\Phi^\lambda(t)=e^{i\lambda_1 t^{\alpha} + i \lambda_2 t^{\beta}} \psi(t)/t,
\eeq
and set $\Phi^\lambda(-t)=0$.
For $a >0$, let
\beq\label{Phi_dfn2}
\Phi_a^\lambda(t)= a^{-1} \Phi^\lambda(t/a).
\eeq 
Let $|\lam| = |\lam_1| + |\lam_2|$.
Then there exists $\gamma_0>0$ such that for all $r \geq 1$ and all $F \in L^2(\R)$,
\[\left\|\sup_{a>0, |\lambda|\ge r} |F*\Phi_a^\lambda| \right\|_{L^2(dx)} \lesssim r^{-\gamma_0} \|F\|_{L^2(dx)}.\]
\end{lemma}
\begin{remark}
We note that as a byproduct of the proof of Lemma \ref{lemma:main}, $\gamma_0$ can be chosen to be independent of $\alpha,\beta$.
\end{remark}

\begin{remark}For $\alpha,\beta\in\N$ this is merely a special case of Stein and Wainger's Theorem 1 in \cite{SteinWainger01}, but to prove Lemma \ref{lemma:main} in full generality requires estimates of a very different flavor. See also work of the first author \cite{Guo15} for a similar result regarding a phase comprised of a single fractional monomial. Fractional exponents appear naturally during the analysis of the operators $B^{m,d}$ via a change of variables $t^m\to t$ (for instance, see \eqref{F_need} and \eqref{eqn:asymmosc}).
(In addition, Theorems \ref{thm:main1} and \ref{thm:main2} could be somewhat generalized to non-integral $m,d$, but we have chosen the integer setting for our main results, to avoid unnecessary complications.) 

\end{remark}

The key contrast of our proof of Lemma~\ref{lemma:main} with the corresponding result in Stein and Wainger \cite{SteinWainger01} appears in the proof of Lemma \ref{lemma:kernelest}. 
The strategy is to linearize the operator $F \mapsto \sup_{a>0, |\lambda|\ge r} |F*\Phi_a^\lambda|$ using stopping-times for $a,\lam$, and to bound an oscillatory integral by showing that for all but a small exceptional region of the integral, the phase has a large derivative of some order.  Our proof enables us to make the exceptional region independent of the precise stopping-time $\lam$, thus obviating the need for the small-set maximal functions that appear in \cite{SteinWainger01}; at the cost of restricting our attention to phases with only two monomials, we are also able to handle fractional powers.

\subsection{Method of proof: Theorem \ref{thm:main2}}\label{sec_method_thm2}
Next, we sketch the proof of Theorem~\ref{thm:main2}. To analyze $A^{m,1}_N$, where $m \geq 3$ is an integer, we first decompose the operator as
\[
A^{m,1}_N = \sum_{k \in \mathbb{Z}} A^{m,1}_N \circ P_k,
\]
where $P_k$ is a Littlewood-Paley projection onto frequency $\sim 2^k$ in the $y$-variable. In view of the modulation invariance (\ref{2809ee1.14}) in the $x$-variable, this is the only viable Littlewood-Paley decomposition we can use for the operator $A^{m,1}$; a Littlewood-Paley decomposition in the $x$-variable is doomed to fail. We also note that the Littlewood-Paley projection in the $y$-variable commutes with $A^{m,1}_N$, since the stopping time $N$ in the operator $A^{m,1}_N$ depends only on $x$ but not on $y$. 

Now to analyze each Littlewood-Paley piece of $A^{m,1}_N$, we decompose the integral \[A^{m,1}_N \circ P_k f(x,y) = \int_{\mathbb{R}} (P_k f)(x-t,y-t^{m}) e^{iN(x) t} \frac{dt}{t}\] 
into two parts, where $t$ is small or large compared to the frequency $2^k$. For $t$ small, we compare the resulting integral to a maximally truncated Carleson operator in the $x$-variable; this is natural in view of the remarks in Section \ref{sec_consequences}. The error will be given by a strong maximal function, since $P_k f$ is localized in frequency in the $y$-variable. For $t$ large, we need to use a van der Corput estimate: again we take advantage of the fact that the stopping time $N$ of $A^{m,1}_N$ depends only on $x$, to take a partial Fourier transform in the $y$-variable.

In order to reassemble the various Littlewood-Paley pieces, the main ingredient is a vector-valued estimate for the maximally truncated Carleson operator (Theorem \ref{lemma:vvtrunccarleson}).

A similar strategy works for $B^{m,m}_N$ for $m \geq 2$ an integer. There is, however, an interesting distinction depending on whether $m$ is odd or even: when $m$ is odd, we need to use the maximally truncated Carleson operator in the $y$-variable, whereas when $m$ is even, the component of the operator that would correspond to the maximally truncated Carleson operator magically vanishes. (See equation \eqref{2709ee3.5}, and the discussion immediately thereafter.)

Roughly speaking, our proof of Theorem \ref{thm:main2} works because the linearizing function depends on the same variable in which the modulation invariance occurs, so the other variable is at our disposal to use Plancherel's theorem and localize in frequency via Littlewood-Paley decomposition. Essential parts of this proof fail in the remaining cases $A^{2,1}$, $A^{m,m}$ and $B^{m,1}$, where $m>1$. 
For $A^{m,m}$, the linearizing function varies with $x$, so we would like to use Plancherel's theorem in $y$ and localize in the $y$ frequency. However, the modulation invariance in \eqref{2809ee1.15} causes translation invariance in the $y$ frequency so that any attempt at doing a Littlewood-Paley decomposition is doomed from the start. Similar behavior occurs for $B^{m, 1}$.

\section{Preliminaries}
\subsection{Notation}
The notation $A\lesssim B$ always means $A\le C\cdot B$ with $0<C<\infty$ depending only on $m,d$ and the function $\psi$ chosen below (and within the proof of Lemma \ref{lemma:main}, on $\alpha,\beta$). Similarly, $A\approx B$ means $C_1 A\le B\le C_2 A$ with $0<C_1\le C_2<\infty$ and the same dependence.
We use the Fourier transform $\hat{f}(\xi) = \int_{\R} f(x) e^{-i \xi x} dx$ with inverse $\check{g} (x) =(2\pi)^{-1} \int_\R g(\xi) e^{i x \xi} d\xi$ and Plancherel identity $\|f\|_2 = (2\pi)^{-1/2} \|\hat{f} \|_2$.

\subsection{Littlewood-Paley decomposition}\label{sec_notation}
Once and for all we fix a smooth function \(\psi:\R\to\R\) supported on \(\{t:1/2~\le~|t|~\le~2\}\) such that $0\le \psi(t)\le 1$ and \(\sum_{k\in\Z} \psi_k(t)=1\) for all \(t\not=0\), where \(\psi_k(t)=\psi(2^{-k} t)\).
Define the associated Littlewood-Paley projection of a function $F$ on $\R$ by
\begin{align}\label{eqn:lwpaley_0}
F_k(w) = P_k F(w)=\RInt F(u) \check{\psi}_k(w-u) du,
\end{align}
where 
$\check{\psi}_k$ denotes the inverse Fourier transform of the function $\psi_k$. 
The standard Littlewood-Paley estimates apply, in the form 
\[ \|F\|_p \lesssim \Big\| \left(\sum_{k} |P_k F|^2 \right)^{1/2} \Big\|_{p} \lesssim \|F\|_{p}.\]
We will apply this in the $x$-variable or $y$-variable of $f(x,y)$, depending which is free.

\subsection{Vector-valued inequalities}
In this section we collect several vector-valued estimates that will play important roles in our work. 

Define the maximally truncated Carleson operator by
\beq\label{define_Cstar}
 \mathscr{C}^* F(x) = \sup_{N\in\R,\varepsilon>0} \Big| p.v.\int_{|t|\le\varepsilon} F(x-t) e^{iNt} \frac{dt}{t}\Big|. 
 \eeq
Note that this operator is usually studied with the inequality $|t|\le\varepsilon$ being reversed; we may of course reduce to that case by subtracting the Carleson operator from $\mathscr{C}^*$.

\begin{letterthm}\label{lemma:vvtrunccarleson}
For $1<p<\infty,$ 
\begin{align}\label{eqn:vvtrunccarleson}
\Big\|\Big(\sum_{k\in\Z} |\mathscr{C}^* F_k|^2\Big)^{1/2}\Big\|_p \lesssim \Big\|\Big(\sum_{k\in\Z} |F_k|^2\Big)^{1/2}\Big\|_p,
\end{align}
with a constant depending only on $p$.
\end{letterthm}
We assemble the necessary results to verify Theorem \ref{lemma:vvtrunccarleson} in a brief appendix (Section \ref{sec:vector_valued_0}).

Next, let $\mathscr{M}$ be the maximal operator of Radon-type along the curve $(t,t^m)$:
\begin{align}\label{eqn:defparmaxfct}
\mathscr{M} f(x,y)=\sup_{r>0} \frac1{2r} \Int_{-r}^{r} |f(x-t,y-t^m)| dt.
\end{align}
This is known to be a bounded operator of $L^p$ for $1<p\leq \infty$ (e.g. by a small modification of the proof in the case of the parabola $(t,t^2)$, \cite[Chapter XI \S1.2, \S 2]{Stein}). We require a vector-valued inequality for $f_k := P_k f$, with $P_k$ acting on either the $x$-variable or $y$-variable (to be specified later):

\begin{letterthm}\label{lemma:vectorvalued}
For $1<p<\infty$ we have
\begin{align}\label{eqn:vvmaxfct}
\Big\|\Big(\sum_{k\in\Z} |\mathscr{M} f_k|^2\Big)^{1/2}\Big\|_p \lesssim \Big\|\Big(\sum_{k\in\Z} |f_k|^2\Big)^{1/2}\Big\|_p,
\end{align}
with a constant depending only on $p$.
\end{letterthm}
This result is stated in \cite[Theorem 2.5]{RdFRT86}, as a consequence obtainable from a more general theory.
For completeness, we offer a brief, self-contained  proof for our special case in an appendix (Section
\ref{sec:vector_valued}); we thank E. M. Stein for sharing with us this method of proof, which appears in a significantly more general form in the preprint \cite[Appendix A, Theorem A.1]{MST15x}.

Finally, we will need two one-variable Hardy-Littlewood maximal functions in the plane, denoted by $M_1$ and $M_2$. Indeed, they will act on the first and second variable respectively:
\begin{align}
M_1 f(x,y)=\sup_{r>0} \frac1{2r} \int_{-r}^r |f(x-u,y)| du \label{eqn:onevarmaxfct1} \\
M_2 f(x,y)=\sup_{r>0} \frac1{2r} \int_{-r}^r |f(x,y-t)| dt \label{eqn:onevarmaxfct2}. 
\end{align}
They are bounded on $L^p(\mathbb{R}^2)$ for all $1 < p < \infty$, and satisfy the following vector-valued inequality, which follows easily by integrating a corresponding result of Fefferman and Stein:
\begin{letterthm}\label{lemma:vvsmaxfct}
For $1<p<\infty$ and $i = 1,2$, we have
\begin{align}
\Big\|\Big(\sum_{k\in\Z} |M_i f_k|^2\Big)^{1/2}\Big\|_p \lesssim \Big\|\Big(\sum_{k\in\Z} |f_k|^2\Big)^{1/2}\Big\|_p,
\end{align}
with a constant depending only on $p$.
\end{letterthm}
See e.g. \cite[Chapter II \S1.1]{Stein} for further details.

\section{The asymmetric case: Theorem \ref{thm:main1}}\label{sec:reduction}\label{sec:asymm}
First we prove Theorem \ref{thm:main1}, assuming Lemma \ref{lemma:main}; then in Section \ref{sec_proof_lemma_main} we prove the lemma.

For convenience, we define the auxiliary variable $z=z(x,y)$ to be understood as indicating either $z(x,y)=x$ or $z(x,y)=y$, so that $N(z)$ can mean either $N(x)$ or $N(y)$. To simplify notations, we also define 
\beq
Tf(x,y)=H^{m,d}_{N(z)} f(x,y),
\endeq
with $m,d$ satisfying the conditions of Theorem \ref{thm:main1}.
With $\psi_\ell$ as defined in Section \ref{sec_notation}, define for each $\ell \in \Z$
\[T_\ell f(x,y)=\int_{\R} f(x-t,y-t^m) e^{iN(z)t^d} \psi_\ell(t) \frac{dt}t.\]
Let \(n:\R\to\Z\) be such that for all $z \in \R$,
\beq\label{eqn:defnx}
2^{-n(z)d}\le |N(z)|< 2^{-(n(z)-1)d}.
\eeq
Then we decompose \(T=T^{(1)}+T^{(2)}\) with 
\[T^{(1)}f(x,y)=\sum\limits_{\ell\le n(z)} T_\ell f(x,y)\] and 
\[T^{(2)}f(x,y)=\sum\limits_{\ell>0} T_{n(z)+\ell} f(x,y).\]
The motivation for this decomposition is that when $\ell \le n(z)$, $\psi_\ell(t)$ localizes to $|t| \leq 2^{\ell+1} \leq 2^{n(z)+1}$ and the exponential factor $e^{iN(z)t^d}$ is well approximated by $1$. Consequently we write $T^{(1)}f(x,y)$ as
\beq\label{T1_two_terms}
\sum_{\ell\le n(z)} \RInt f(x-t,y-t^m) (e^{iN(z)t^d}-1) \psi_\ell(t) \frac{dt}t + \sum_{\ell\le n(z)}\RInt f(x-t,y-t^m) \psi_\ell(t) \frac{dt}t.
\eeq
We may estimate the absolute value of the first summand brutally by applying \eqref{eqn:defnx}:
\[\lesssim \sum_{\ell\le n(z)} \RInt |f(x-t,y-t^m)|\cdot |N(z)t^{d-1} \psi_\ell(t)| dt\lesssim \frac{1}{2^{n(z)+2}} \Int_{-2^{n(z)+1}}^{2^{n(z)+1}} |f(x-t,y-t^m)| dt.\]
The right hand side is bounded by $\mathscr{M}f(x,y)$, where $\mathscr{M}$ denotes the  maximal operator along $(t,t^m)$ defined in (\ref{eqn:defparmaxfct}).

The second summand in (\ref{T1_two_terms}) is bounded in absolute value by the maximal truncated Hilbert transform along the curve $(t,t^m)$, defined by
\begin{align}\label{eqn:defmaxht}
\mathcal{H}^* f(x,y)=\sup_{\varepsilon, R > 0}\Big|\int_{\varepsilon < |t| < R} f(x-t,y-t^m) \frac{dt}t\Big|,
\end{align}
plus an error term bounded by $\mathscr{M}f(x,y)$ (which arises at the endpoint when passing from smooth bump functions to a sharp truncation).
Thus in total we have obtained the pointwise estimate 
\[|T^{(1)}f|\lesssim \mathscr{M}f + \mathcal{H}^* f.\] 
Since both $\mathcal{H}^*$, $\mathscr{M}$ are known to be bounded in $L^p$, $1<p<\infty$ (for example, by slight modifications of Stein and Wainger's work for $(t,t^2)$ in \cite{SteinWainger78}),  we may conclude that
 \[\|T^{(1)} f\|_p\lesssim \|f\|_p\] for all $1<p<\infty$.

 It remains to show the same for $T^{(2)}$.
Let \[S_\ell f(x,y)=T_{n(z)+\ell} f(x,y);\] 
we claim that it suffices to prove there exists some $\ga_0>0$ such that for all $\ell > 0$,
\begin{align}\label{eqn:l2est}
\|S_\ell f\|_2\lesssim 2^{-\gamma_0 \ell} \|f\|_2.
\end{align}
Indeed, the triangle inequality implies the pointwise estimate $|S_\ell f|\lesssim \mathscr{M}f$, so that we immediately obtain $\|S_\ell f\|_p\lesssim \|f\|_p$ for all $1<p<\infty$; by interpolation  with (\ref{eqn:l2est}) we then obtain for any $1<p<\infty$ there exists some $\ga_p>0$ such that
\[\|S_\ell f\|_p\lesssim 2^{-\gamma_p \ell} \|f\|_p.\]
 Finally, summing over $\ell \geq 0$ gives
  \[\|T^{(2)}f\|_p\lesssim \|f\|_p.\]

All that remains is to prove \eqref{eqn:l2est}; we proceed by distinguishing two cases.

\subsection{The $B_N^{m,d}$ operators.}\label{sec:asymmBcase}
Here we consider the case $z(x,y)=y$. Applying Plancherel's theorem in the free $x$-variable, we obtain
\[\|S_{\ell} f(x,y)\|_{L^2(dx)} = (2\pi)^{-1/2} \Big\| \RInt g_\xi(y-t^m) e^{iN(y)t^d - i\xi t} \psi_{n(y)+\ell}(t) \frac{dt}t\Big\|_{L^2(d\xi)}, \]
where \[g_\xi(y)=\int_\R e^{-i\xi x} f(x,y) dx.\]
Therefore to prove (\ref{eqn:l2est}) it will suffice to prove a bound of the form 
\[ \Big\| \RInt F(y-t^m) e^{iN(y)t^d-i\xi t} \psi_{n(y)+\ell}(t)\frac{dt}t \Big\|_{L^2(dy)} \lesssim 2^{-\gamma_0 \ell} \|F\|_2,  \]
uniformly in \(\xi\in\R\), for all single variable functions $F$. Recall that the cutoff function $\psi_{n(y)+\ell}$ has supports both in the positive half line and in the negative half line. Accordingly let us split the integration over $t$ into a positive and a negative part. 
We consider the positive part; the negative component is treated in an entirely analogous way.
Changing variables $t^m \mapsto t$, we see that it suffices to show there exists some $\ga_0>0$ such that for all $\ell>0$ and all $F \in L^2(\R)$,
\beq\label{F_need}
 \Big\|\int\limits_{0}^\infty F(y-t) e^{iN(y)t^{d/m}-i\xi t^{1/m}} \psi_{n(y)+\ell}(t^{1/m}) \frac{dt}t\Big\|_{L^2(dy)}\lesssim 2^{-\gamma_0 \ell} \|F\|_2, 
 \eeq
uniformly in $\xi$.
In fact (\ref{F_need}) is an immediate consequence of the key Lemma \ref{lemma:main}, with $\alpha=d/m$, $\beta=1/m$.
To see this, we first re-write $N(y) = 2^{-n(y)d +r(y)d}$ with $0<r(y)<1$ for all $y$. 
Then 
for $a \in \R$ and $\lam \in \mathbb{R}^2$, we define $\Phi_a^\lam := a^{-1}\Phi^\lam(t/a)$, where
\[
\Phi^\lam(t) := e^{i \lam_1 t^{\alpha} + i \lam_2 t^{\beta}} \psi(t^{1/m}) t^{-1} \quad \text{for $t > 0$},
\]
and $\Phi^\lam(t) = 0$ for $t \leq 0$.
One then observes that the integral on the left hand side of (\ref{F_need}) is equal to $F* \Phi_{a}^\lam(y)$, with parameters 
\[ a=2^{(n(y) + \ell)m}, \qquad \lam_1 = 2^{\ell d + r(y) d}, \qquad \lam_2=  -\xi 2^{n(y)+\ell}.\]
 Then (recalling $\ell >0$, $0<r(y)<1$), we have
\[ |\lam| = |\lam_1| + |\lam_2| \geq 2^{\ell d + r(y) d} \geq 2^{\ell d}, \] and we see from Lemma \ref{lemma:main} that for any fixed $\ell >0$, 
\[ \Big\|\int\limits_{0}^\infty F(y-t) e^{iN(y)t^{d/m}-i\xi t^{1/m}} \psi_{n(y)+\ell}(t^{1/m}) \frac{dt}t\Big\|_{L^2(dy)}\lesssim \| \sup_{a>0, |\lam| \geq 2^{\ell d}} |F * \Phi_a^\lam | \|_2 \lesssim 2^{-\gamma_0 \ell} \|F\|_2, \]
as desired. This proves \eqref{F_need} and hence \eqref{eqn:l2est} in this case.

\subsection{The $A_N^{m,d}$ operators.}\label{sec:asymmAcase} Here we treat the case $z(x,y)=x$.
Applying Plancherel's theorem in the  free $y$-variable, we obtain
\[\Big\|S_\ell f(x,y)\Big\|_{L^2(dy)}=(2\pi)^{-1/2}\Big\|\int\limits_{\R} g_\eta(x-t) e^{iN(x)t^d-i\eta t^m} \psi_{n(x)+\ell}(t) \frac{dt}t\Big\|_{L^2(d\eta)}\]
with \[g_\eta(x)=\int_\R e^{-i\eta y} f(x,y) dy.\]
By Plancherel's theorem it suffices to show that there exists $\ga_0>0$ such that for each $\ell>0$, 
$$\Big\|\int\limits_{\R} F(x-t) e^{iN(x)t^d-i\eta t^m} \psi_{n(x)+\ell}(t) \frac{dt}t\Big\|_{L^2(dx)}\lesssim 2^{-\gamma_0 \ell} \|F\|_2$$
uniformly in $\eta$. 
Now it is clear that we may proceed similarly to (\ref{F_need}), and deduce this bound from 
Lemma \ref{lemma:main} with $\alpha=d$, $\beta=m$.

\subsection{Proof of Lemma \ref{lemma:main}}\label{sec_proof_lemma_main}
 In order to complete the proof of Theorem~\ref{thm:main1}, it remains to prove Lemma \ref{lemma:main}. Due to a minor technical issue we will assume the pair $\{\alpha,\beta\}\not=\{2,3\}$ in the proof. However, this case is of course already covered by Stein and Wainger's work \cite[Theorem 1]{SteinWainger01}.

In fact it suffices to prove there exists $\ga_0$ such that for all $r \geq 1$,
\begin{align}\label{eqn:penult}
\|\sup_{a>0, \; r\le |\lambda|\le 2r} |F*\Phi_a^\lambda| \|_2 \lesssim r^{-\gamma_0} \|F\|_2.
\end{align}
For with this result in hand, we immediately obtain the desired result,
\[ \Big\| \sup_{a>0, \; |\lambda|\ge r} |F*\Phi_a^\lambda| \Big\|_2 \le \sum_{k=0}^\infty  \Big\| \sup_{a>0,\; 2^{k}r \le |\lambda|\le 2^{k+1} r } |F*\Phi_a^\lambda| \Big\|_2 \lesssim r^{-\gamma_0} \|F\|_2. \]

We proceed by linearizing the supremum. For measurable functions $a:\R\to(0,\infty)$, $\lambda:\R\to\R^2$ with $r\le |\lambda(u)|\le 2r$ for all $u\in\R$, we define an operator $\Lambda:L^2(\R)\to L^2(\R)$ by
\[ \Lambda F(u)=F*\Phi^{\lambda(u)}_{a(u)}(u)=\int_\R F(t) \Phi_{a(u)}^{\lambda(u)} (u-t) dt. \]
The bound \eqref{eqn:penult}  will follow from proving $\|\Lambda\|_{2\to 2}\lesssim r^{-\gamma_0}$
 for some $\gamma_0>0$ with the implicit constant independent of $a,\lambda$. Since $\|\Lambda\|_{2\to 2}=\|\Lambda \Lambda^*\|_{2\to 2}^{1/2},$ we will in fact prove
 \beq\label{Lambda22}
  \|\Lambda \Lambda^*\|_{2\to 2}\lesssim r^{-2\gamma_0}.
  \eeq
We calculate
\begin{align}
\Lambda \Lambda^* F(u) 
&= \int_\R F(s) (\Phi^\nu_{a_1}*\tilde{\Phi}^\mu_{a_2}) (u-s) ds,\label{eqn:TTstar}
\end{align}
with $\tilde{\Phi}(u):=\overline{\Phi(-u)}$ and $\nu=\lambda(u), \mu=\lambda(s), a_1=a(u), a_2=a(s)$.
Note that by rescaling we may write
\[ (\Phi^\nu_{a_1} * \tilde{\Phi}^\mu_{a_2})(s)= a_2^{-1} (\Phi^\nu_{a_1/a_2}*\tilde{\Phi}^\mu_1)(a_2^{-1}s)=a_1^{-1} (\Phi^\nu_{1}*\tilde{\Phi}^\mu_{a_2/a_1})(a_1^{-1}s).\]
Thus we will deduce (\ref{Lambda22}) from applying the following bounds, which are the heart of the proof:
\begin{lemma}\label{lemma:kernelest}There exists $\gamma_1>0$ such that for $0< h\le 1$, $r\le |\nu|, |\lambda|\le 2r$ we have
\begin{align}
\label{eqn:kernelest}
|(\Phi^\nu_h * \tilde{\Phi}^\mu_1)(s)| \lesssim r^{-\gamma_1} \mathbf{1}_{\{|s|\le 4\}}(s) + \mathbf{1}_{\{|s|\le r^{-\gamma_1}\}}(s), \\
\label{eqn:kernelestB}
|(\Phi^\nu_1 * \tilde{\Phi}^\mu_h)(s)| \lesssim r^{-\gamma_1} \mathbf{1}_{\{|s|\le 4\}}(s) + \mathbf{1}_{\{|s|\le r^{-\gamma_1}\}}(s).
\end{align}
\end{lemma}
\begin{remark}Note that the exceptional sets in \eqref{eqn:kernelest}, \eqref{eqn:kernelestB} do not depend on $\nu,\mu$. This is in contrast to \cite[Lemma 4.1]{SteinWainger01}. As a consequence we do not require Stein and Wainger's small set maximal function \cite[Proposition 3.1]{SteinWainger01}.\end{remark}

We first proceed with the proof of Lemma \ref{lemma:main}, and then prove Lemma \ref{lemma:kernelest} in Section \ref{sec_lemma_main_sub}.
Applying \eqref{eqn:kernelest} and \eqref{eqn:kernelestB} appropriately (depending on whether $a_1 \geq a_2$ or $a_1 \leq a_2$), we deduce
\[
|(\Phi^\nu_{a_1} * \tilde{\Phi}^\mu_{a_2})(s)| \lesssim r^{-\gamma_1} \sum_{k=1}^2 \left(a_k^{-1} \mathbf{1}_{\{|s|\le 4 a_k\}}(s) + (a_k r^{-\gamma_1})^{-1}\mathbf{1}_{\{|s|\le r^{-\gamma_1} a_k\}}(s) 
\right).\]
Thus for any $G \in L^2$ we may compute
\begin{multline*}
|\langle \Lambda \Lambda^* F,G\rangle| = \Big| \int_\R\int_\R (\Phi^\nu_{a_1} * \tilde{\Phi}^\mu_{a_2})(u-s) F(s) \overline{G(u)} ds du\Big|
	\\ \lesssim 
 r^{-\gamma_1} \left( \int_\R MF(u) |G(u)| du + \int_\R |F(s)| MG(s) ds \right),
\end{multline*}
in which $M$ denotes the standard one-variable Hardy-Littlewood maximal function. 
(Here the important point is that we have integrated first in whichever variable was independent of the stopping-time $a_k$, for the two terms $k=1,2$.)
Via the Cauchy-Schwarz inequality and the boundedness of $M$ on $L^2$, we obtain
\[ |\langle \Lambda \Lambda^* F,G\rangle|\lesssim r^{-\gamma_1} \|F\|_2 \|G\|_2. \]
This completes the proof of \eqref{eqn:penult} with $\gamma_0=\gamma_1/2$. \\

\subsection{Proof of Lemma \ref{lemma:kernelest}}\label{sec_lemma_main_sub}
 We will only prove \eqref{eqn:kernelest}, as (\ref{eqn:kernelestB}) follows by symmetry.
 By definition,
\begin{align}\label{2709ee2.40}
(\Phi^\nu_h * \tilde{\Phi}^\mu_1)(s) 
&= \int_\R  e^{i\nu_1 t^\alpha+i\nu_2 t^\beta-i\mu_1(ht-s)^\alpha-i\mu_2 (ht-s)^\beta} \frac{\psi(t)}t  \frac{\overline{\psi(ht-s)}}{ht-s} dt.
\end{align}
First notice that the support of $\Phi^\nu_h*\tilde{\Phi}^\mu_1 (s)$ is contained in $\{s:|s|\le 4\}$.

In order to apply van der Corput estimates, we need to analyze when the phase function
\[ Q(t,s)=\nu_1 t^\alpha + \nu_2 t^\beta - \mu_1(ht-s)^\alpha - \mu_2(ht-s)^\beta \]
has a large derivative of some order. Here we recall that $\nu_i=\lam_i(u)$ are fixed with respect to $t,s$ (the relevant variables of integration in \eqref{eqn:TTstar}), and that $ r \leq  |\nu_1| + |\nu_2| \leq 2 r$. On the other hand, $\mu_i = \lam_i(s)$ depends on $s$ (in an unknown way), and thus our strategy is to make our argument independent of $\mu_1,\mu_2$.

{\underline{Case 1:} Suppose that $0<h\le h_0$, where $0<h_0<1$ is to be determined later, depending on $r, \al, \be$; this is the easier case.

Let $0<\varepsilon_1<1$ be small and fixed. 
Within the support of $\psi(t)\overline{\psi(ht-s)}$, we estimate
\[ |\partial_t Q(t,s)| \ge |\alpha \nu_1 t^{\alpha-1} + \beta\nu_2 t^{\beta-1} | - h C r, \]
where $C$ is a positive constant only depending on the exponents $\alpha,\beta$.
Let us define the function
\[F(t)=\alpha \nu_1 t^{\alpha-1}+\beta\nu_2 t^{\beta-1},\]
and its associated exceptional set
\[ E=\{t\in[1/2,2]\,:\,|F(t)|\le \tau r^{1-\varepsilon_1}\},\] where $\tau$ is a positive constant that depends only on $\alpha,\beta$ and is to be determined later. Our strategy will be to choose $\tau$ so that $|E|$ is small and then apply van der Corput's lemma outside of $E$. 

We will prove (at the end of the considerations for Case 1):
\begin{lemma}\label{lemma_E}
There exists a choice of $\tau$ (depending only on $\al,\be$) such that 
\begin{equation}\label{2709ee2.45}
|E|\lesssim r^{-\varepsilon_2}
\end{equation}
with $\ep_2 = \ep_1/|\be-\al|$.
\end{lemma}

Assuming $\tau$ is chosen as in the lemma, we now specify $h_0$ to be such that 
\begin{equation}\label{2709ee2.48}
h_0 C = \frac12 \tau r^{-\varepsilon_1}.
\end{equation} 
Then whenever $h \leq h_0$,  for all $t \in [1/2,2] \setminus E$,
\beq\label{2709ee2.49}
|\partial_t Q(t,s)| \gtrsim r^{1-\varepsilon_1}.
\endeq 
We now split the integral in \eqref{2709ee2.40} according to whether $t\in [1/2,2]\setminus E$ or $t\in E$. 
We estimate the portion of the integral over $E$ trivially by the measure of $E$, which is small $\lesssim r^{-\ep_2}$ by Lemma \ref{lemma_E}.

We will estimate the portion of the integral over $[1/2,2] \setminus E$ by applying van der Corput's lemma  combined with the lower bound \eqref{2709ee2.49}. 

Here we encounter a delicate point: as stated in \cite[Chapter VIII \S 1.2]{Stein} van der Corput's lemma for a first derivative assumes monotonicity. We circumvent this assumption as follows. We first note that $E$ (and thus also $[1/2,2]\backslash E$) is a finite union of intervals, with the number of intervals being bounded by a small absolute constant. 
To see this note that the equation
\[ \al \nu_1 t^{\al-1} + \be \nu_2 t^{\be-1} \pm \tau r^{1-\ep_1} = 0\]
has at most $3$ solutions in $t>0$ (see for example \cite[Lemma 3]{SteinWainger70}).

Thus we may apply the following slight variant of van der Corput's lemma (proved at the end of Case 1) to each such interval:
\begin{lemma}\label{lemma_vdC}
Suppose $\phi$ is real-valued and smooth in $(a,b)$ and that both $|\phi' (x)| \geq \sig_1 $ and $|\phi'' (x) | \leq \sig_2$ for all $t \in(a,b)$. Then 
\[\left| \int_a^b e^{i \lam \phi(t)}dt \right| \leq (a-b) \left( \frac{\sig_2}{\sig_1^2}\right) \lam^{-1}.\]
\end{lemma}
Here we note that for $s$ fixed, we have that $Q(t,s)$ is $C^\infty$ with respect to $t$ for all $t$ in the support of $\psi(t)\overline{\psi(ht-s)}$; in particular note that both $t, ht-s$ are bounded away from the origin. We also verify trivially that for all such $t$, 
\beq\label{upperQ}
 | \partial_t^2 Q(t,s)|  \lesssim r, 
 \eeq
with a constant depending only on $\al,\be$. 
Hence applying Lemma \ref{lemma_vdC} with the bounds (\ref{2709ee2.49}) and (\ref{upperQ}) to each of the finitely many finite-length intervals in $[1/2,2] \setminus E$, we obtain for each such portion of the integral a bound of size $\lesssim r (r^{1-\ep_1})^{-2} = r^{-(1-2\ep_1)}$. 
In total, combining this with our trivial estimate for the portion of the integral over $E$, we have proved
\[|(\Phi_h^\nu*\tilde{\Phi}^\mu_1)(s)| \lesssim r^{-(1-2\varepsilon_1)} + r^{-\varepsilon_2}\lesssim r^{-\varepsilon_3},\]
for all $|s| \leq 4$, for a suitable $\ep_3 >0$, which suffices for (\ref{eqn:kernelest}) in this case.
All that remains is to verify Lemmas \ref{lemma_E} and \ref{lemma_vdC}.

\begin{proof}[Proof of Lemma \ref{lemma_E}]
We observe that if one of $|\nu_1|,|\nu_2|$ dominates the other then  $|F(t)|$ is large, that is $|F(t)|\gtrsim r$. 
More precisely, recall that $|\nu_1| + |\nu_2| \approx r$, and suppose that $ |\nu_2|/|\nu_1| \leq c_0$ for some small constant $c_0$ (so in particular $|\nu_1| \gtrsim r$). Then 
\[ |F(t)|= |\nu_1| \left| \al t^{\al-1} + \be \frac{\nu_2}{\nu_1} t^{\be-1} \right|;\]
if $c_0$ is chosen sufficiently small (with respect to $\al,\be$) we may guarantee that for all $t \in [1/2,2]$,
\[ \al t^{\al-1} \geq 2 \left| \be \frac{\nu_2}{\nu_1} t^{\be-1} \right| \]
and hence 
\[ |F(t)| \geq  |\nu_1|  \frac{\al}{2} t^{\al-1} \geq c_1 r,\]
say. 
We may argue similarly to obtain $|F(t)| \geq c_1'r$ if $|\nu_1|/|\nu_2| \leq c_0'$ where $c_0'$ depends only on $\al,\be$. 

By choosing $\tau < \min\{c_0,c_0'\}$ (hence depending only on $\al,\be$) we then see that if $E$ is to be non-empty, we must be in the regime where 
$c_0^{-1} \leq |\nu_1|/|\nu_2| \leq c_0'$, that is, $|\nu_1| \approx |\nu_2|$.
In this case, we will deduce that (\ref{2709ee2.45}) holds.
Suppose that $\al< \be$;  write $c:= \nu_2/\nu_1$ so that $|c| \in [c_0^{-1},c_0'].$
Then
\[ F(t)=\alpha \nu_1 t^{\alpha-1} (1+ c(\be/\al) t^{\beta-\alpha}),\] 
 so that for all $t\in E$ we must have
\[ r|1+c (\be/\al)t^{\beta-\alpha}| \leq |F(t)| \leq \tau r^{1-\varepsilon_1},\] 
that is, $t$ must satisfy
\[ |1+c (\be/\al)t^{\beta-\alpha}| \leq \tau r^{-\ep_1}.\]
The measure of such $t$ is $\lesssim r^{-\ep_1/(\be-\al)}$, with an implicit constant dependent on $\al,\be$. 
For the case $\alpha>\beta$ we argue in an entirely analogous way. This proves Lemma \ref{lemma_E}.
\end{proof}
\begin{proof}[Proof of Lemma \ref{lemma_vdC}] We recall the proof of the original van der Corput lemma in the case of a first derivative \cite[Ch VIII, Proposition 2]{Stein}, which bounds the integral in question by 
\[ \lam^{-1} \int_a^b \left| \frac{d}{dt} \left( \frac{1}{\phi'} \right) \right| dt  = \lam^{-1} \int_a^b \left| \frac{\phi''(t)}{\phi'(t)^2} \right| dt,\]
where we have evaluated the derivative rather than invoking monotonicity of $\phi'$ to bring the absolute values outside the integral. The inequality claimed in Lemma \ref{lemma_vdC} then clearly follows. 
\end{proof}

We have now concluded the proof of Lemma \ref{lemma:kernelest} in Case 1.

{\underline{Case 2.}} In the remaining case, $h_0 \leq h \leq 1$. Fix any small $0<\ep_4<1$; if $|s|\le r^{-\ep_4}$ we use the triangle inequality to bound \eqref{2709ee2.40} trivially by 1, which is sufficient for the second term in (\ref{eqn:kernelest}).
Thus from now on we assume that 
\beq\label{s_big}
  |s|\ge r^{-\ep_4}
  \eeq
and work to obtain a small bound for the integral.
Note that as a vector,
\begin{align}\label{matrix_equation}
\left(\begin{matrix} \partial_t Q(t,s)\\
\partial_t^2 Q(t,s)\\
\partial_t^3 Q(t,s)\\
\partial_t^4 Q(t,s)\end{matrix}\right)
 = M_{t,s} \left(\begin{matrix}
\alpha\nu_1 t^{\alpha-1}\\
-\alpha\mu_1 (ht-s)^{\alpha-1}\\
\beta\nu_2 t^{\beta-1}\\
-\beta\mu_2 (ht-s)^{\beta-1}
\end{matrix}\right),
\end{align}
where $M_{t,s}$ is the $4 \times 4$ matrix
\begin{align}\label{eqn:matrix}
M_{t,s} = \left(\begin{matrix}
1 & h & 1 & h\\
a_1 t^{-1} & a_1 h^2(ht-s)^{-1} & b_1 t^{-1} & b_1 h^2(ht-s)^{-1}\\
a_2 t^{-2} & a_2 h^3(ht-s)^{-2} & b_2 t^{-2} & b_2 h^3(ht-s)^{-2}\\
a_3 t^{-3} & a_3 h^4(ht-s)^{-3} & b_3 t^{-3} & b_3 h^4(ht-s)^{-3}
\end{matrix} \right), 
\end{align}
\noindent and
\begin{align*}
a_1=\alpha-1, \quad a_2&=(\alpha-1)(\alpha-2), \quad a_3=(\alpha-1)(\alpha-2)(\alpha-3),\\
b_1=\beta-1, \quad b_2&=(\beta-1)(\beta-2), \quad b_3=(\beta-1)(\beta-2)(\beta-3).
\end{align*}

 If we can show that $|\det M_{t,s}|$ is sufficiently large, that is 
 \beq\label{M_big}
| \det M_{t,s} | \gtrsim r^{-\kappa}
\eeq
 for some $\kappa>0$, then we will apply the following lemma (whose proof we defer to the end of the section):
\begin{lemma}\label{lemma:detest}
Let $A$ be an invertible $n\times n$ matrix and $x\in\R^n$. Then
\[ |Ax| \ge |\det A| \|A\|^{1-n} |x|, \]
where $\|A\|$ denotes the matrix norm $\sup_{|x|=1} |Ax|$.
\end{lemma}

Note that $\|M_{t,s} \| \lesssim 1$ (since we only consider $t$ in the support of $\psi(t)\overline{\psi(ht-s)}$, so that both $t, ht-s$ are bounded away from the origin).  If we have shown (\ref{M_big}) for $t$ in a certain interval, then applying Lemma \ref{lemma:detest} to (\ref{matrix_equation}), we see that throughout that interval,
\beq\label{sum_Q}
\left( \sum_{k=1}^4 |\partial_t^k Q(t,s)|^2 \right)^{1/2} \gtrsim r^{-\kappa} |(\nu_1,\nu_2,\mu_1,\mu_2)^T| \gtrsim r^{1-\kappa}.
 \eeq

Then applying the van der Corput lemma to that portion of the integral (\ref{2709ee2.40}) shows that portion is bounded by $r^{-(1-\kappa)/4}$. (Note: to be precise, if only the first order term $|\partial_t Q(t,s)|$ dominates in (\ref{sum_Q}), then we must apply the variant Lemma \ref{lemma_vdC} of the van der Corput lemma, using the trivial upper bound $|\partial_t^2Q(t,s)| \lesssim r$, similar to our argument in Case 1. This will result in a bound for the portion of the integral over that interval of size $\lesssim r^{-(1-2\kappa)}$, which is sufficient.)

In fact, we will show that $|\det M_{t,s}|$ is sufficiently large in this manner for all but a small exceptional set $E$ of $t$, with measure $ \lesssim r^{-\kappa'}$ for some small $\kappa'>0$. (As in our argument in Case 1, we will also note that this exceptional set is a union of a finite number of intervals, dependent only on $\al,\be$, so that we may apply the above argument to each individual component of $[1/2,2]\setminus E$.) Thus this strategy is sufficient to complete the proof of (\ref{eqn:kernelest}).

We require the following purely algebraic identity.
\begin{lemma}\label{lemma:det}
Let $a_0,\dots,a_3,b_0,\dots,b_3,x,y$ be arbitrary real numbers. Then
\[ \left|\begin{matrix} 
a_0 & a_0 & b_0 & b_0\\
a_1 x & a_1 y & b_1 x & b_1 y\\
a_2 x^2 & a_2 y^2 & b_2 x^2 & b_2 y^2\\
a_3 x^3 & a_3 y^3 & b_3 x^3 & b_3 y^3
\end{matrix}\right| = (c(x^2+y^2)+dxy)(x-y)^2 xy,\]
where $c,d$ are given by
\begin{align}\label{eqn:detconst}
c = -\left|\begin{matrix}a_0 & b_0\\a_1 & b_1\end{matrix}\right|\cdot
\left|\begin{matrix}a_2 & b_2\\a_3 & b_3\end{matrix}\right|, \qquad
d = c+ \left|\begin{matrix}a_0 & b_0\\a_3 & b_3\end{matrix}\right|\cdot
\left|\begin{matrix}a_1 & b_1\\a_2 & b_2\end{matrix}\right|.
\end{align}
\end{lemma}

\begin{proof}
Expand the determinant along the first row, combine the terms corresponding to the first and second columns, and those corresponding to the third and fourth columns, respectively.
Then again expand each of the two resulting $3\times 3$ determinants along the first row.
\end{proof}

Let $\tilde{M}$ denote the matrix in Lemma \ref{lemma:det}. Rescaling the individual rows and columns of $M_{t,s}$ appropriately to clear denominators, we see that 
\[ \det M_{t,s}= h^2 t^{-6}(ht-s)^{-6} \det  \tilde{M}\]
where within $\tilde{M}$ we set $a_0=b_0=1,$ $x=ht-s$ and $y=ht$. Then we may apply 
Lemma \ref{lemma:det} to compute
\begin{align}\label{eqn:detM}
\det M_{t,s} =t^{-5} (ht-s)^{-5} s^2 h^3 S(t),
\end{align}
with 
\[ S(t)=h^2(2c+d)t^2 - h(2c+d)st + cs^2,\]
in which $c,d$ are as in (\ref{eqn:detconst}).
Note that with $a_0=b_0=1$ and the other $a_i,b_i$ as specified above, then $c\not=0$ is equivalent to $\alpha\not=\beta$, $\alpha,\beta\not=1,2$. 

In order to now verify that $|\det M_{t,s}|$ is sufficiently large, as in (\ref{M_big}), we  distinguish between two cases.

{\underline{Case 2A.}} Suppose first that $2c+d=0$, so that $S(t)=cs^2$. We must then verify that $c \neq 0$. 
Since $2c+d=0$, clearly if $c=0$ then $d=0$. But recall from above that $c = 0$  implies that either $\al=2$ or $\be=2$ (since the hypotheses of Lemma \ref{lemma:main} already ruled out $\al=\be$, $\al=1$ or $\be=1$). Recall also that we assume in this stage of the proof that the pair $(\al,\be)$ is not $(2,3)$.

Suppose that $\al=2$. Then we would have $d = (\beta-1)^2(\beta-2)^2(\beta-3)(\alpha-1),$ which is clearly non-zero (since $\be \neq 3$). Analogously we see that $\be=2$ leads to a contradiction. Thus we may conclude that $c\neq 0$, and recalling $|s| \geq r^{-\ep_4}$ from (\ref{s_big}) and  $h \geq h_0 \gtrsim r^{-\ep_1}$ from \eqref{2709ee2.48}, we may compute immediately from  (\ref{eqn:detM}) that 
 \[|\det M_{t,s} |\gtrsim r^{-4\ep_4-3\ep_1},\]
holds for all $t \in [1/2,2]$. This verifies (\ref{M_big}) and allows us to apply the van der Corput lemma to bound the full integral (\ref{2709ee2.40}) by $r^{-\kappa}$ for some $\kappa>0$, completing the proof of (\ref{eqn:kernelest}) in this case.

\underline{Case 2B.} The final case we must consider is when $2c+d\not=0$.
Fix any small $0<\ep_5<1$ and define 
\[ E=\{t\in[1/2,2]\,:\, |S(t)|\le r^{-2\ep_5-2\varepsilon_1}\}.\]
Note first that $E$ is a union of at most two intervals, since $S$ is a quadratic polynomial.
Then for $t \in [1/2,2] \setminus E$, 
\eqref{eqn:detM} in combination with $|s| \geq r^{-\ep_4},$ $h \gtrsim r^{-\ep_1}$ implies
\[|\det M_{t,s} | \gtrsim r^{-2\ep_4-5\ep_1-2\ep_5},\]
verifying (\ref{M_big}) so that we may apply the van der Corput lemma to bound the portion of the integral over $[1/2,2] \setminus E$ by $\lesssim r^{-\kappa}$ for some $\kappa>0$, which suffices for the first term in (\ref{eqn:kernelest}).

We will bound the portion of the integral over $E$ trivially, so all that remains is to verify that $E$ has small measure, for which we call upon the following lemma (see Christ \cite[Lemma 3.3]{Christ85}):
\begin{lemma}\label{lemma:sublevelset}
Let $I\subset\R$ be an interval, $k\in\N$, $f\in C^k(I)$, and suppose that for some $\sig>0$, $|f^{(k)}(x)|\ge \sig$ for all $x\in I$. Then there exists a constant $0<C<\infty$ depending only on $k$ such that for every $\rho >0$,
\begin{align}\label{eqn:sublevelset}
|\{x\in I\,:\, |f(x)|\le \rho\}|\le C \Big(\frac\rho \sig \Big)^{1/k}.
\end{align}
\end{lemma}
By the choice of $h_0$ in \eqref{2709ee2.48} we have 
\[|S''(t)|\gtrsim h^2 \gtrsim h_0^2 \gtrsim r^{-2\varepsilon_1}.\] 
Thus by Lemma \ref{lemma:sublevelset}, we have $|E|\lesssim r^{-\ep_5}$, which suffices for the second term in (\ref{eqn:kernelest}).

All that remains to complete the proof of Lemma \ref{lemma:kernelest}, and hence of the main Lemma \ref{lemma:main}, is to verify Lemma \ref{lemma:detest}.
\begin{proof}[Proof of Lemma \ref{lemma:detest}]
First we show $\|A^{-1}\|\le \|A\|^{n-1}/|\det A |$. By homogeneity we can assume $\|A\|=1$. Then all the eigenvalues of $AA^*$ are between $0$ and $1$. Let $\lambda$ be the smallest eigenvalue of $AA^*$. Then $\|A^{-1}\|=\lambda^{-1/2}\le \det(AA^*)^{-1/2}=|\det A|^{-1}.$
Therefore in general,
\[ |x|=|A^{-1}Ax|\le \|A^{-1}\|\cdot  |Ax| \le \|A\|^{n-1} |\det A |^{-1} |Ax|,\]
as desired.
\end{proof}

\section{The symmetric case: Theorem \ref{thm:main2}}\label{sec:symm}

Here we prove Theorem \ref{thm:main2}. We present  the proof in detail only for $B^{m,m}$; thus in the following we write $T = B^{m,m}$. The proof for $A^{m,1}$ is, mutatis mutandi, analogous, and we merely sketch the necessary changes in Section \ref{sec:remarksproofA}.

Recalling the function $\psi_k$ fixed in Section \ref{sec_notation}, we define the Littlewood-Paley projection in the free $x$-variable by
\begin{align}\label{eqn:lwpaley_1}
f_k(x,y) = P_k f(x,y)=\RInt f(u,y) \check{\psi}_k(x-u) du,
\end{align}
where $\check{\psi}_k$ denotes the inverse Fourier transform of the function $\psi_k$.
In particular, note that $T P_k=P_k T$.

\subsection{Single annulus estimate.}
We fix $k_0 \in \Z$ and split the operator as $T=T_{k_0}^{(1)}+T_{k_0}^{(2)}$, where for any fixed $k,$ $T^{(1)}_k$ is defined as
\[ T_k^{(1)}f(x,y) := p.v. \int_{|t|\le 2^{-k}} f(x-t,y-t^m) e^{iN(y)t^m} \frac{dt}t \]
and $T_k^{(2)} := T - T_{k}^{(1)}$ accordingly. 

Our  key estimate for $T^{(1)}_{k_0}$ is the pointwise bound:
\beq\label{T1_ineq}
	|T^{(1)}_{k_0} P_{k_0} f(x,y)|\lesssim \mathscr{C}^*P_{k_0} f(x,y) + M_1 M_2 P_{k_0} f(x,y),
	\eeq
in which the maximally truncated one-variable Carleson operator $\mathscr{C}^*$ is defined as in (\ref{define_Cstar}); here our understanding is that $\mathscr{C}^*$ acts only on the second variable. Also, $M_1$ and $M_2$ refer to the Hardy-Littlewood maximal function in the first and second variable respectively, as defined in \eqref{eqn:onevarmaxfct1}, \eqref{eqn:onevarmaxfct2}.

We will prove (\ref{T1_ineq}) by using the fact that under the  single  annulus assumption 
\beq\label{single_annulus}
f=P_{k_0} f,
\eeq
the function $f(x-t,y-t^m)$ is well approximated by $f(x,y-t^m)$.
Precisely, we assume (\ref{single_annulus}) and estimate
 \[|T_{k_0}^{(1)}f(x,y)|\le \mathrm{\bf I}+\mathrm{\bf II},\] where
\begin{eqnarray}
\mathrm{\bf I}&=&\int\limits_{|t|\le 2^{-k_0}} |f(x-t,y-t^m)-f(x,y-t^m)| \frac{dt}{|t|}, \label{eqn:symlowerfreq} \\
\mathrm{\bf II}&=&\Big| p.v.\int\limits_{|t|\le 2^{-k_0}} f(x,y-t^m) e^{iN(y)t^m} \frac{dt}{t}\Big|.\label{2709ee3.5}
\end{eqnarray} 
At this point there is a striking dichotomy in our treatment, depending on the parity of $m$: if $m$ is even, the term $\mathrm{\bf II}$ vanishes identically due to the integrand being an odd function. On the other hand, if $m$ is odd, we can change variables $t^m\mapsto t$ (appropriately in the cases $t>0, t<0$) to see
\begin{align}\label{eqn:singleannulcarleson}
\mathrm{\bf II}\lesssim \sup_{\varepsilon>0}\Big| p.v. \int\limits_{|t|\le \varepsilon} f(x,y-t) e^{iN(y)t} \frac{dt}{t}\Big| 
	\leq \mathscr{C}^* f(x,y),
\end{align}
in which the maximally truncated Carleson operator acts only on the second variable. This contributes the first term in (\ref{T1_ineq}).

Next, we note that the first term $\mathrm{\bf I}$ can be estimated by a maximal function due to the single annulus assumption (\ref{single_annulus}). 
We write
\begin{align}\label{eqn:singleannulustailcalc}
f(x-t,y-t^m)-f(x,y-t^m)=\int_{\R} f(x-u,y-t^m)(\check{\psi}_{k_0}(u-t)-\check{\psi}_{k_0}(u)) du.
\end{align}
By the rapid decay of the first derivative of $\check{\psi}$ we certainly have 
\beq\label{decay_xi}
|\frac{d}{d\xi}(\psi_{k_0})\check{\;} | = |\frac{d}{d\xi} (2^{k_0} \check{\psi}(2^{k_0}\xi) )|  \leq 2^{2k_0} (1+|2^{k_0}\xi|)^{-2}.
\eeq
Now suppose that for some $j \geq 0$, $u$ is in the annulus
 \beq\label{u_annulus}
 2^{-k_0+j-1}\le |u|\le 2^{-k_0+j},\eeq
so that for $|t| \leq 2^{-k_0}$ we have both $2^{-k_0+j-2}\le |u|, |u-t|\le 2^{-k_0+j+1}.$
Thus applying the mean value theorem and the decay (\ref{decay_xi}), for $u$ in the annulus (\ref{u_annulus}) we have 
\[ |\check{\psi}_{k_0}(u-t)-\check{\psi}_{k_0}(u)|\lesssim |t|  \cdot 2^{2(k_0-j)},\]
where the implicit constant depends only on the choice of $\psi$.

Therefore \eqref{eqn:singleannulustailcalc} can be estimated in absolute value by 
\beq\label{eqn:singleannuluscalc1} \lesssim |t| 2^{2k_0} \sum_{j=0}^\infty 2^{-2j} \int_{|u|\le 2^{-k_0+j}} |f(x-u,y-t^m)| du. \eeq
This allows us to bound the term $\mathrm{\bf I}$ by 
\beq\label{2809ee3.11}
\lesssim \sum_{j=0}^\infty 2^{-j} \frac{1}{2^{-k_0} 2^{-k_0+j}} \int_{|u|\le 2^{-k_0+j}}\int_{|t|\le 2^{-k_0}}  |f(x-u,y-t^m)| dt du .
\eeq
We may dominate this by the maximal functions $M_1$ and $M_2$ as follows. 
Indeed, we focus temporarily on the inner integration in $t$ in (\ref{2809ee3.11}):
\begin{align*}
& \frac1{2^{-k_0}}\int_{|t|\le 2^{-k_0}}  |f(x-u,y-t^m)| dt
\leq & \frac{C}{2^{-k_0}}\int_{|s|\le 2^{-k_0 m}} |f(x-u,y-s)| |s|^{\frac1m -1} ds
\end{align*}
Since $|s|^{\frac{1}{m} - 1} \mathbf{1}_{|s|\le 2^{-k_0 m}}$ is radially decreasing and integrable in $s$, with integral equal to $C 2^{-k_0}$, by \cite[Chapter III Theorem 2]{Stein_Sing},
we can bound the above by $\lesssim M_2f(x-u,y)$.
Thus \eqref{2809ee3.11} is bounded by $M_1 M_2 f(x,y)$. This completes the proof of the inequality (\ref{T1_ineq}) for $T^{(1)}_{k_0}$.

We now turn to estimating $T_{k_0}^{(2)}f$, still under the single annulus assumption (\ref{single_annulus}). 
 Let us define for any integer $\ell$
\[ T_{\ell} f(x,y)=\int_{\R} f(x-t,y-t^m) e^{iN(y)t^m} \psi_{\ell}(t)\frac{dt}{t}. \]
Then certainly
\beq\label{2809ee3.14}
|T_{k_0}^{(2)}f|\lesssim \mathscr{M}f + \sum_{\ell=0}^\infty |T_{-k_0+\ell} f|;\eeq
here we need merely observe  that the maximal operator $\mathscr{M}$ along $(t,t^m)$ (defined in (\ref{eqn:defparmaxfct})) arises in (\ref{2809ee3.14}) due to the transition to smooth cutoffs.

Next, we claim that (still under the assumption (\ref{single_annulus})) there exists a constant $\ga >0$ such that for all $\ell \geq 0,$
\begin{align}\label{eqn:symmcase-l2decay}
 \|T_{-k_0+\ell} f\|_2 \lesssim 2^{-\gamma \ell} \|f\|_2.
\end{align}
To prove \eqref{eqn:symmcase-l2decay} we proceed similarly to Section \ref{sec:asymmBcase}. First we apply Plancherel's theorem in the free $x$-variable, 
so that it is equivalent to prove that for $g_\xi(y)=\int_\R e^{-i \xi x} f(x,y) dx,$
\[
\Big\|\int\limits_{\R} g_\xi(y-t^m) e^{iN(y)t^m-i\xi t} \psi_{-k_0+\ell}(t) \frac{dt}t\Big\|_{L^2 (d\xi, dy)}\lesssim 2^{-\gamma \ell} \|g_\xi\|_{L^2(d\xi, dy)}.
\]
In particular, we note that due to the assumption (\ref{single_annulus}), $g_\xi$ is nonzero only in the frequency annulus $2^{k_0-1}\le |\xi|\le 2^{k_0+1}$.

We then split the integral into a positive and negative part, which are dealt with analogously. We focus here on the positive portion of the integral; by a change of variables $t^m\mapsto t$ the claim (\ref{eqn:symmcase-l2decay}) is reduced to showing
\begin{align}\label{eqn:asymmosc}
\Big\|\int\limits_{0}^\infty F(y-t) e^{iN(y)t-i\xi t^{1/m}} \psi_{-k_0+\ell}(t^{1/m}) \frac{dt}t\Big\|_{L^2(dy)}\lesssim 2^{-\gamma \ell} \|F\|_{L^2(dy)}
\end{align}
for all single variable functions $F$, uniformly in 
$2^{k_0-1}\le |\xi|\le 2^{k_0+1}$. 

As in the proof of Lemma \ref{lemma:main} we proceed by the $TT^*$ method. 
For convenience we write \[\tilde{\psi}_k(t)=\psi_k(t^{1/m}) \mathbf{1}_{(0,\infty)}(t),\] 
and denote the operator on the left hand side of \eqref{eqn:asymmosc} by $\tilde{T}$.
Then $\|\tilde{T}\|_{2\to 2}=\|\tilde{T}\tilde{T}^*\|_{2\to 2}^{1/2}$, where
\beq\label{TTK}
\tilde{T}\tilde{T}^* F(y)=\int_\R F(y-s) K_{N(y),N(y-s)}(s) ds
\eeq
and for any $\lam_1,\lam_2\in \R$ the kernel $K_{\lam_1,\lam_2}$ is given by
\[ K_{\lambda_1,\lambda_2}(s)=\int_\R e^{i\lambda_1 t-i\lambda_2 (t-s)-i\xi (t^{1/m}- (t-s)^{1/m})} \frac{\tilde{\psi}_{-k_0+\ell}(t-s)}{t-s} \frac{\tilde{\psi}_{-k_0+\ell}(t)}{t} dt. \]
Via the substitution $t\mapsto \rho t$ with $\rho=2^{m(-k_0+\ell)}$ we obtain
\[ \rho K_{\lambda_1,\lambda_2}(\rho s)=\int_\R e^{i\lambda_1\rho t-i\lambda_2\rho(t-s)-i\xi2^{-k_0+\ell}(t^{1/m}- (t-s)^{1/m})} \frac{\tilde{\psi}_{0}(t-s)}{t-s} \frac{\tilde{\psi}_{0}(t)}{t} dt. \]
We need to analyze the phase function
\begin{align}\label{eqn:symmphasefct}
Q(t,s)=\lambda_1\rho t-\lambda_2\rho (t-s) + \eta (t^{1/m}-(t-s)^{1/m}),
\end{align}
where $\eta=-\xi 2^{-k_0+\ell}$, so in particular
\beq\label{eta_bound}
2^{\ell-1}\le |\eta|\le 2^{\ell+1}.
\eeq
 On first sight this may not look promising, because the phase function includes linear terms which tend to cause trouble (compare Stein and Wainger \cite{SteinWainger01}). However, it turns out that we are allowed to take derivatives to isolate the non-linear term (recall $m \geq 2$) because we know by (\ref{eta_bound}) that its coefficient $\eta$ is large. Taking two derivatives with respect to $t,$ we obtain
\[ \partial^2_t Q(t,s)=c\eta (t^{\alpha}-(t-s)^{\alpha}) \]
with $\alpha=\frac1m-2$ and $c=\frac1m (\frac1m-1)$.
Suppose that $|s|\ge 2^{-\ell/2}$. Then by the mean value theorem  and (\ref{eta_bound}), $|\partial^2_t Q(t,s)|\gtrsim 2^{\ell/2}$ throughout the region of  $t$ and $t-s$ considered (depending only on the support of $\tilde{\psi}_0$).
In this case, an application of the second derivative test shows that
\[|\rho K_{\lambda_1,\lambda_2}(\rho s)|\lesssim 2^{-\ell/4}.\] 
On the other hand, if $|s|\le 2^{-\ell/2}$ we merely use the triangle inequality for the trivial estimate
\[|\rho K_{\lambda_1,\lambda_2}(\rho s)|\lesssim 1. \]
Altogether we have proved, with $\rho=2^{m(-k_0+\ell)}$,
\[|K_{\lambda_1,\lambda_2}(s)|\lesssim 2^{-\ell/4} \rho^{-1} \mathbf{1}_{\{|s|\le 4\rho\}}(s) + \rho^{-1}\mathbf{1}_{\{|s|\le 2^{-\ell/2}\rho\}} (s), \]
uniformly in $\lambda_1,\lambda_2$. 
Applying this in (\ref{TTK}) allows us to deduce that
\beq\label{2709ee3.26}
 |\tilde{T}\tilde{T}^* F(y)| \lesssim 2^{-\ell/4} MF(y), \eeq
where $MF$ denotes the standard one-variable Hardy-Littlewood maximal function.
An application of the $L^2$ estimate for $M$ now implies our claim \eqref{eqn:asymmosc} with $\gamma=1/8$;
by Plancherel we then finally obtain (\ref{eqn:symmcase-l2decay}).

\subsection{Square function estimate}

In this section we assemble the single annulus estimates of the previous section to derive the $L^p$ boundedness of our operator $T=B^{m,m}$. This  application of the Littlewood-Paley theory is in the spirit of Bateman and Thiele \cite{BatemanThiele}.

In view of the relation $TP_k=P_kT$ and the standard Littlewood-Paley inequalities, we have
\[ \|Tf\|_p \lesssim \Big\|\left(\sum_{k\in\Z} |TP_k f|^2\right)^{1/2}\Big\|_p\lesssim \Big\|\left(\sum_{k\in\Z} |T^{(1)}_k P_k f|^2\right)^{1/2}\Big\|_p+\Big\|\left(\sum_{k\in\Z} |T^{(2)}_k P_k f|^2\right)^{1/2}\Big\|_p\]

In the term for $T^{(1)}_k$ on the right hand side  we apply the estimate (\ref{T1_ineq}); then by applying the vector-valued estimates of 
Theorems   \ref{lemma:vvtrunccarleson}  and  \ref{lemma:vvsmaxfct} to the maximally truncated Carleson operator 
and the one-variable maximal function,  we  obtain
\[ \Big\|\left(\sum_{k\in\Z} |T^{(1)}_k P_k f|^2\right)^{1/2}\Big\|_p\lesssim \Big\|\left(\sum_{k\in\Z} |P_k f|^2\right)^{1/2}\Big\|_p\lesssim \|f\|_p. \]

For $T_k^{(2)}$, we recall that by \eqref{2809ee3.14}
\[|T^{(2)}_k P_k f|\lesssim \mathscr{M}P_k f + \sum_{\ell=0}^\infty |T_{-k+\ell} P_k f|,\]
so that by Minkowski's inequality for integrals,
\[  \Big\|\left(\sum_{k\in\Z} |T^{(2)}_k P_k f|^2\right)^{1/2}\Big\|_p 
	\lesssim
	 \Big\|\left(\sum_{k\in\Z} | \mathscr{M} P_k f|^2\right)^{1/2}\Big\|_p 
	 	+ 
		\sum_{\ell=0}^\infty \Big\|\left(\sum_{k\in\Z} |T_{-k+\ell} P_k f|^2\right)^{1/2}\Big\|_p .
	\]
We may apply Theorem \ref{lemma:vectorvalued} to obtain a bound $\lesssim \|f\|_p$  for the first term; for the second term it would suffice to show that for each $1<p<\infty$ there exists $\ga>0$ such that for every $\ell \geq 0,$
\begin{align}\label{eqn:symmpenult}
\Big\|\left(\sum_{k\in\Z} |T_{-k+\ell} P_k f|^2\right)^{1/2}\Big\|_p\lesssim 2^{-\gamma\ell}\Big\|\left(\sum_{k\in\Z} |P_k f|^2\right)^{1/2}\Big\|_p
	\lesssim 2^{-\gamma \ell} \|f\|_p.
\end{align}
 For $p=2$ this follows from \eqref{eqn:symmcase-l2decay}. But recall that we always have the simple estimate $|T_{\ell} P_k f|\lesssim \mathscr{M} P_k f$ for all $k,\ell\in\Z$ by the triangle inequality. Therefore Theorem \ref{lemma:vectorvalued} implies a bound without decay, namely
\beq\label{eqn:symmpenult'}
\Big\|\left(\sum_{k\in\Z} |T_{-k+\ell} P_k f|^2\right)^{1/2}\Big\|_r\lesssim \|f\|_r,
\eeq
valid for all $1 < r < \infty$. 
Now in general, \eqref{eqn:symmpenult} follows for all $1< p < \infty$ by interpolating between the $L^2 \to L^2(\ell^2)$ case of (\ref{eqn:symmpenult}) and the $L^r \to L^r(\ell^2)$ bound (\ref{eqn:symmpenult'}) for the vector-valued map $f \mapsto \{T_{-k+\ell} P_k f\}_{k \in \mathbb{Z}}$. The proof of Theorem \ref{thm:main2} is now complete.

\subsection{Remarks on the proof for $A^{m,1}$}\label{sec:remarksproofA}
As mentioned above, we will not repeat the proof explicitly for $A^{m,1}$, but merely complement the sketch already provided in Section \ref{sec_method_thm2} by pointing out two key modifications. Of course one interchanges the roles of the $x$ and $y$ variables. In addition:
\begin{enumerate}
\item[(1)] The cancellation miracle for even $m$ in the term $\mathrm{\bf II}$ of (\ref{2709ee3.5}) does not occur for $A^{m,1}$. Instead one always needs to invoke Carleson's theorem in the form of Theorem \ref{lemma:vvtrunccarleson}, analogous to the computation already carried out for odd $m$ in (\ref{eqn:singleannulcarleson}).
\item[(2)] In the treatment of $A^{m,1}$, the restriction $m\not=2$ originates because the relevant phase function analogous to (\ref{eqn:symmphasefct}) is
\[ 
Q(t,s) = \lam_1 \rho t - \lam_2 \rho (t-s) + \eta(t^m - (t-s)^m).
\]
Visibly, when $m=2$, the phase function $Q(t, s)$ is now \emph{linear} in $t$, so that its second derivative vanishes, and consequently we fail in this case to obtain a good bound for the kernel. 
\end{enumerate}

\section{Deductions for partial Carleson operators}\label{sec_deductions}
\subsection{$L^2$ consequences of partial Carleson bounds}\label{sec_deductions1}

As stated in Section \ref{sec_consequences}, the $L^2$ boundedness of $A^{m,1}_N$, for any fixed integer $m \geq 1$, implies the $L^2$ boundedness of Carleson's operator (\ref{eq:carleson}). Similarly, the $L^2$ boundedness of $B^{m,m}_N$, when $m \geq 1$ is an odd integer, implies the $L^2$ boundedness of Carleson's operator.  (Of course,  in the work of this paper, our logic is actually the other way round: in proving Theorem~\ref{thm:main2}, we used Carleson's theorem as a black box.)

We will see how to carry out these deductions from a more general argument we now give in the context of
 the quadratic Carleson operator $\mathscr{C}^{\text{par}}$ along the parabola (defined in equation (\ref{eqn:paraboliccarleson})).
 We prove that an inequality of the form (\ref{partial_HPN}) would imply the analogue over $\R$ of Lie's result \cite{Lie09} on the one-variable quadratic Carleson operator $\mathscr{C}_Q$:

 \begin{prop} \label{prop_CQ}
 Assume the veracity of the estimate
\begin{align}
\label{eqn:partialparabolicL2}
\Big\|\Int_{\epsilon\le |t|\le R} f(x-t,y-t^2)  e^{iN_1(x)t+iN_2(x)t^2}\frac{dt}{t}\Big\|_{L^2(dxdy)}\le C\|f\|_2,
\end{align}
for all Schwartz functions $f$, where $N_1,N_2: \R \maps \R$ are measurable functions, $0<\epsilon<R$ are real parameters, and the constant $C$ is independent of $f,N_1,N_2,\epsilon,R$. 
Then the operator 
\[ f \mapsto \mathscr{C}_Qf(x) := \sup_{N\in\R^2} \Big|p.v.\int_\R f(x-t) e^{iN_1t+iN_2t^2} \frac{dt}t\Big|\]
is bounded on $L^2(\R)$. 
\end{prop}
Note that in our assumed bound (\ref{eqn:partialparabolicL2}), the  linearizing  functions $N_1,N_2$ are independent of $y$, so this is a far weaker assumption than the conjectured $L^2$ bound for $\mathscr{C}^{\text{par}}$ in (\ref{eqn:paraboliccarleson}).
In the argument that we will now give for (\ref{eqn:partialparabolicL2}), if we replace the curve $(t,t^2)$  by $(t,t^m)$ and  the phase by $N_1(x)t + N_2(x)t^m$, and furthermore specify that $N_2$ is identically zero, we may deduce Carleson's original theorem from the partial bound for $A^{m,1}_N$ for any integer $m \geq 1$; or, if we specify $N_1$ is identically zero, we may deduce Carleson's original theorem from the partial bound for $B^{m,m}_N$ for $m$ an odd integer. (When $m$ is even, under the specification $N_1 \con 0$, the operator in (\ref{eqn:partialparabolicL2}) would vanish, due to the integrand being an odd function.)

In general, to prove Proposition \ref{prop_CQ}, we use an elementary tensor 
$f(x,y)=h(x)g(y)$, where $h,g$ are real Schwartz functions, in which case (\ref{eqn:partialparabolicL2}) implies
\[ \Big\|\Int_{\epsilon\le |t|\le R} h(x-t)e^{iN_1(x)t+iN_2(x)t^2} g(y-t^2) \frac{dt}{t}\Big\|_{L^2(dx dy)}\le C\|h\|_2 \|g\|_2. \]
Applying Plancherel's theorem in the $y$ variable we obtain
\beq\label{eqn:quadraticcarleson1}
 \Big\|\Int_{\epsilon\le |t|\le R} h(x-t)e^{iN_1(x)t+iN_2(x)t^2} \widehat{g}(\eta)e^{-i\eta t^2} \frac{dt}{t}\Big\|_{L^2(dx d\eta)}\le  C\|h\|_2\|g\|_2. \eeq
Suppose for the time being that we have chosen $g$ such that we have an estimate of the form
\beq\label{eqn:quadraticcarlesonerr} \Big\|\Int_{\epsilon\le |t|\le R} h(x-t)e^{iN_1(x)t+iN_2(x)t^2} \widehat{g}(\eta)(e^{-i\eta t^2}-1)\frac{dt}{t}\Big\|_{L^2(dx d\eta)}\le C \|h\|_2 \|g\|_2. 
\eeq
We would deduce from (\ref{eqn:quadraticcarleson1}) and (\ref{eqn:quadraticcarlesonerr}) that
\[ \Big\|\int_{\epsilon\le |t|\le R} h(x-t)e^{iN_1(x)t+iN_2(x)t^2} \widehat{g}(\eta) \frac{dt}t\Big\|_{L^2(dx d\eta)}\le C\|h\|_2\|g\|_2, \]
so that by Plancherel and Fubini,
\[ \Big\|\int_{\epsilon\le |t|\le R} h(x-t)e^{iN_1(x)t+iN_2(x)t^2} \frac{dt}t\Big\|_2\le C\|h\|_2. \]
Via Fatou's lemma this gives the $L^2$ boundedness of the quadratic Carleson operator $h \mapsto \mathscr{C}_Q h$, as claimed in Proposition \ref{prop_CQ}.

To obtain the estimate \eqref{eqn:quadraticcarlesonerr}, we choose $\delta$ with $0<\delta<1/R^2$ and specify that $g$ be a Schwartz function on $\R$ such that $\widehat{g}$ is supported on $[-\delta,\delta]$ and $\|g\|_2>0$.
Then by Minkowski's inequality and Fubini, the left hand side of \eqref{eqn:quadraticcarlesonerr} is bounded by
\begin{align}\label{eqn:eqn20150615-1}
\|h\|_2\Int_{\epsilon\le |t|\le R} \|\widehat{g}(\eta)(e^{-i\eta t^2}-1)\|_{L^2(d\eta)} \frac{dt}{|t|}.
\end{align}
 
The mean value theorem, followed by Plancherel, shows that for $|t| \leq R$,
$$\Int_{\R}\Big|\widehat{g}(\eta)(e^{-i\eta t^2}-1)\Big|^2 d\eta\le \delta^2 t^4 \|g\|_2^2.$$
This implies that \eqref{eqn:eqn20150615-1} is no greater than
\beq
\|h\|_2 \|g\|_2 \cdot \delta \Int_{\epsilon\le |t|\le R} |t| dt\le \|h\|_2\|g\|_2,
\endeq
which completes the proof of (\ref{eqn:quadraticcarlesonerr}), and hence Proposition \ref{prop_CQ}.

\subsection{$L^2$ deductions for partial Carleson operators}\label{sec_deductions2}
Remark \ref{remark_difficult} stated that $L^2$ bounds for $A_N^{2,1}$, $A^{m,m}_N$ and $B^{m,1}_N$ (with $m>1$) follow from known Carleson theorems. We briefly indicate these deductions, which are along the lines of arguments in Sections \ref{sec:asymmBcase} and \ref{sec:asymmAcase}.
By Plancherel's theorem in the free $y$-variable,
\[ \|A_N^{m,m} f \|_{L^2 (dxdy)} = \Big\| \int_{\R} g_\eta (x-t) e^{i (N(x) -\eta) t^m} \frac{dt}t \Big\|_{L^2(dx d\eta)}\]
where $g_\eta(x)=\int_\R e^{-i\eta y} f(x,y) dy.$
Then an $L^2$ bound of the form
\[ \Big\| \int_{\R} g_\eta (x-t) e^{i (N(x) -\eta) t^m} \frac{dt}t \Big\|_{L^2(dx)} \lesssim \| g_\eta \|_{L^2(dx)},\]
uniform in $\eta$, follows from Stein and Wainger \cite{SteinWainger01} (since $m>1$), and this suffices.

In the next case,
\[ \|A_N^{2,1} f \|_{L^2 (dxdy)} = \Big\| \int_{\R} g_\eta (x-t) e^{i N(x) t - i \eta t^2} \frac{dt}t \Big\|_{L^2(dx d\eta)}.\] 
Observe that
\[ i\eta t^2 = i\eta (x-t)^2 - i\eta x^2 + 2i\eta x t. \]
Define $Q_\eta f(x)=e^{i\eta x^2} f(x)$ and set $\widetilde{N}(x)=N(x)-2\eta x$.
Then,

\[\int_\R g_\eta(x-t) e^{iN(x)t-i\eta t^2} \frac{dt}{t}=e^{i\eta x^2} \int_\R Q_{-\eta} g_\eta(x-t) e^{i\widetilde{N}(x)t} \frac{dt}{t}=Q_{\eta} H_{\widetilde{N}(x)} Q_{-\eta} g_\eta(x),\]
where $H_N f(x)=\int_\R f(x-t) e^{iNt} \frac{dt}t$. Since $Q_\eta$ is an isometry in $L^2$, our claim follows from the $L^2$ bound for the Carleson operator.

In the final case, by Plancherel's theorem in the free $x$-variable,
\[ \|B_N^{m,1} f \|_{L^2 (dxdy)} = \Big\| \int_{\R} g_\xi (x-t^m) e^{i (N(x) -\eta) t} \frac{dt}t \Big\|_{L^2(d\xi dy)}\]
where
$g_\xi(x)=\int_\R e^{-i\xi x} f(x,y) dx.$
Thus the required $L^2$ bound follows from sending $t \mapsto t^{1/m}$ and applying Guo \cite{Guo15} to the resulting operator, which has one fractional monomial in the phase.

\section{Appendix: Proof of vector-valued inequalities}

\subsection{Proof of Theorem \ref{lemma:vvtrunccarleson}}\label{sec:vector_valued_0}
We assemble results from the Grafakos texts \cite{Grafakos1}, \cite{Grafakos2}. By \cite[Lemma 6.3.2]{Grafakos2}, there is a positive constant $c>0$ such that for any $1\leq p<\infty$, for all $f \in L^p(\R)$ we have the pointwise inequality
\beq\label{CCM}
 \mathscr{C}^* f \leq c M f + M( \mathscr{C} f),
 \eeq
in which $M$ is the standard one-dimensional Hardy-Littlewood maximal function. Since the vector-valued $L^p(\ell^2)$ inequality analogous to (\ref{eqn:vvtrunccarleson}) is known to hold for the Hardy-Littlewood maximal function 
(see e.g. \cite[Chapter II \S1.1]{Stein}), the problem is then reduced to proving the analogue of (\ref{eqn:vvtrunccarleson}) for the Carleson operator $\mathscr{C}$. In fact, this is a special case of   \cite[Exercise 6.3.4]{Grafakos2}, which claims that for all $1<p,r<\infty$ and all weights $w \in A_p$,
\beq\label{Cvv}
\Big\| \left( \sum_k | \mathscr{C} f_k|^r \right)^{1/r} \Big\|_{L^p(w)} \leq C_{p,r}(w) \Big\| \left( \sum_k | f_k|^r \right)^{1/r} \Big\|_{L^p(w)} 
\eeq
for all sequences of functions $f_k \in L^p(w)$. 
This inequality may be verified, following Grafakos, by the method of extrapolation.  We need only note that \cite[Theorem 6.3.3]{Grafakos2} provides a weighted estimate 
\beq\label{C_statement}
\| \mathscr{C}f\|_{L^p(w)} \leq C(p,[w]_{A_p}) \|f \|_{L^p(w)},
\eeq
for every $1<p<\infty$ and $w \in A_p$. 
This is sufficient to prove (\ref{Cvv}) for all the stated values of $r,p$ by applying the vector-valued extrapolation result \cite[Corollary 7.5.7]{Grafakos1}. (Here we remark on the detail that Corollary 7.5.7, to which we appeal, requires that $C(p,[w]_{A_p})$ be an increasing function in $[w]_{A_p}$. We can insure this is the case if we have the  statement, slightly stronger than (\ref{C_statement}), that for every $B>0$ there exists a constant $C_p(B)$ such that for all $w \in A_p$ with $[w]_{A_p} \leq B$ we have $\| \mathscr{C}f \|_{L^p(w)} \leq C_p(B) \|f\|_{L^p(w)}$; such a statement is verified by the explicit version of (\ref{C_statement}) given by Lerner and Di Plinio \cite[Theorem 1.1]{DiPlinioLerner14}.)

Alternatively, once one has the pointwise inequality (\ref{CCM}) and has consequently reduced matters to proving an $L^p(\ell^2)$ vector-valued inequality for $\mathscr{C}$, one can turn to the original result \cite{RdFRT86} in the $L^p(\ell^2)$ case, or the recent streamlined proof \cite[Theorem 7.1]{DemeterSilva15}. 

\subsection{Proof of Theorem \ref{lemma:vectorvalued}}\label{sec:vector_valued}
We recall that the scalar-valued $L^p$-bound for $\mathscr{M}$ was obtained by comparing it to a square function \cite[Chapter XI \S1.2]{Stein}. Indeed, let $\chi(t)$ be a non-negative smooth function with compact support on the interval $[-2,2]$, such that $\chi(t) \equiv 1$ on $[-1,1]$. For $k \in \mathbb{Z}$, let $\chi_k(t) = 2^{-k} \chi(2^{-k}t)$, $d\mu_k(x,y) = \delta_{y=x^m} \chi_k(x)$, and
$$
A_k f(x,y) = f*d\mu_k(x,y) = \int_{\mathbb{R}} f(x-t,y-t^m) \chi_k(t) dt.
$$
Also let $\phi(x,y)$ be a smooth function with compact support on the unit ball in $\mathbb{R}^2$, normalized such that 
$$
\int_{\mathbb{R}^2} \phi(x,y) dxdy = \int_{\mathbb{R}} \chi(t) dt.
$$
For $k \in \mathbb{Z}$, let $\phi_k(x,y) = 2^{-(m+1)k} \phi(2^{-k}x,2^{-mk}y)$, and
$$
B_k f(x,y) = f*\phi_k(x,y) = \int_{\mathbb{R}^2} f(x-u,y-v) \phi_k(u,v) dudv.
$$
Then for non-negative functions $f$, we have the pointwise inequality
\begin{equation} \label{eq:sqdom}
\mathscr{M}f \leq \sup_{k \in \mathbb{Z}} B_k f + Sf,
\end{equation}
where $S$ is the following square function:
\begin{equation} \label{eq:sqmaxdef}
Sf := \left( \sum_{k \in \mathbb{Z}} |A_k f - B_k f|^2 \right)^{1/2}.
\end{equation}
Now $\sup_{k \in \mathbb{Z}} B_k f$ is bounded by the standard maximal function associated to non-isotropic `squares' of sizes $R \times R^m$ on $\mathbb{R}^2$. 
It is known that a vector-valued estimate holds for the maximal function associated to these non-isotropic squares; that is an analogue of Theorem \ref{lemma:vvsmaxfct}. Thus the inequality (\ref{eqn:vvmaxfct}) of Theorem \ref{lemma:vectorvalued} holds for $1 < p < \infty$ if we have $\sup_{k \in \mathbb{Z}} B_k$ in place of $\mathscr{M}$ on the left-hand side. Hence to prove the desired form of (\ref{eqn:vvmaxfct}), all we need to do is to establish
\begin{equation} \label{eq:sqRadonvect} 
\left\| \left( \sum_{\ell \in \mathbb{Z}} |S f_{\ell}|^2 \right)^{1/2} \right\|_{L^p} \lesssim_p  \left\| \left( \sum_{\ell \in \mathbb{Z}} |f_{\ell}|^2 \right)^{1/2} \right\|_{L^p}
\end{equation}
where $S$ is defined by (\ref{eq:sqmaxdef}), and $1 < p < \infty$.

The following scalar-valued inequality for $1<p<\infty$ is already known  \cite[Section 4, Theorem 11]{Stein}:
\begin{equation} \label{eq:sqRadonscalar}
\|Sf\|_{L^p} \lesssim_p \|f\|_{L^p}.
\end{equation}
But to deduce (\ref{eq:sqRadonvect}) we will instead use a related scalar-valued inequality for a signed operator. For $\epsilon_k$  a random sequence of signs $\pm 1$, define
$$
Tf := \sum_{k \in \mathbb{Z}} \epsilon_k (A_k f - B_k f).
$$
It is known that 
\begin{equation} \label{eq:TsqRadonscalar}
\|Tf\|_{L^p} \lesssim_p \|f\|_{L^p}
\end{equation}
for all $1 < p < \infty$, independent of the signs $\epsilon_k$. At the end of this section, we briefly recall a proof of this, for which one uses crucially the non-vanishing of the curvature of the curve $(t,t^m)$, but we first deduce (\ref{eq:sqRadonvect}) from (\ref{eq:TsqRadonscalar}). 

To do so, note that since $T$ is linear, the Marcinkiewicz-Zygmund theorem implies that
$$
\|\, |T f_{\ell}|_{\ell^2} \,\|_{L^p} \lesssim_p \| \, |f_{\ell}|_{\ell^2} \,\|_{L^p}
$$
for $1 < p < \infty$, i.e.
$$
\left\|\, \left|\sum_{k \in \mathbb{Z}} \epsilon_k (A_k f_{\ell}-B_k f_{\ell}) \right|_{\ell^2(d\ell)} \,\right\|_{L^p} \lesssim_p \| \, |f_{\ell}|_{\ell^2} \,\|_{L^p}.
$$
(We write $\ell^2(d\ell)$ to emphasize that the $\ell^2$ norm is taken with respect to the variable $\ell$.) Now we take the expectation, denoted $\mathbb{E}$, over all the possible choices of $\epsilon_k$; by Khintchine's inequality,
$$
 \left( \sum_{k \in \mathbb{Z}} |A_k f_{\ell}-B_k f_{\ell}|^2 \right)^{1/2} \simeq \mathbb{E} \left|\sum_{k \in \mathbb{Z}} \epsilon_k (A_k f_{\ell}-B_k f_{\ell}) \right| .
$$
Taking the $\ell^2(d\ell)$ and then $L^p$ norms on both sides, we get
\begin{align*}
\left\|\, \left( \sum_{k,\ell \in \mathbb{Z}} |A_k f_{\ell}-B_k f_{\ell}|^2 \right)^{1/2} \,\right\|_{L^p} 
&\simeq \left\|\, \left(\mathbb{E} \left|\sum_{k \in \mathbb{Z}} \epsilon_k (A_k f_{\ell}-B_k f_{\ell}) \right|\right)_{\ell^2(d\ell)} \,\right\|_{L^p} \\
&\leq \mathbb{E} \left\|\,  \left|\sum_{k \in \mathbb{Z}} \epsilon_k (A_k f_{\ell}-B_k f_{\ell}) \right|_{\ell^2(d\ell)} \,\right\|_{L^p} \\
&\lesssim_p \mathbb{E} \| \, |f_{\ell}|_{\ell^2} \,\|_{L^p} \\
&= \| \, |f_{\ell}|_{\ell^2} \,\|_{L^p}.
\end{align*}
(The first inequality is the Minkowski inequality.) The left-hand side above is precisely $\|\,| Sf_{\ell} |_{\ell^2} \,\|_{L^p}$. This proves (\ref{eq:sqRadonvect}), and hence (\ref{eqn:vvmaxfct}) of Theorem \ref{lemma:vectorvalued}, for $1 < p <\infty$. 

There are at least two ways of proving (\ref{eq:TsqRadonscalar}).  One is by complex interpolation,  along the lines of arguments in \cite[Section 4]{SteinWainger78},  which we will not discuss here.  Alternatively, we can deduce (\ref{eq:TsqRadonscalar}) from a result of Duoandikoetxea and Rubio de Francia \cite{DuoandikoetxeaRubiodeFrancia86} without using complex interpolation. To do so, let 
$$
d\sigma_k = \epsilon_k ( d\mu_k - \phi_k dxdy ).
$$
Then $d\sigma_k$ has total mass $\|d\sigma_k\| \lesssim 1$, and its Fourier transform satisfies
$$
|\widehat{d\sigma_k}(\xi,\eta)| \lesssim \min\{2^k \|(\xi,\eta)\|,(2^k \|(\xi,\eta)\|)^{-1/m}\}. 
$$
(Here we see the curvature of $(t,t^m)$.)
Furthermore, the operator $\sup_{k \in \mathbb{Z}} |f*|d\sigma_k||$ is bounded by the maximal Radon transform along the curve $(t,t^m)$ plus the Hardy-Littlewood maximal operator adapted to certain non-isotropic balls in $\mathbb{R}^2$. It follows that $\sup_{k \in \mathbb{Z}} |f*|d\sigma_k||$ is bounded on $L^q(\mathbb{R}^2)$ for all $1 < q < \infty$. Thus Theorem B of  Duoandikoetxea and Rubio de Francia \cite{DuoandikoetxeaRubiodeFrancia86} applies, and shows that $Tf = \sum_{k \in \mathbb{Z}} f * d\sigma_k$ is bounded on $L^p$ for all $1 < p < \infty$. This completes our proof of (\ref{eq:TsqRadonscalar}).

\subsection*{Acknowledgements.} We would like to thank C. Thiele and E. M. Stein for many helpful comments and discussions. 
Pierce is supported in part by NSF DMS-1402121. Yung is supported in part by the Hong Kong Research Grant Council Early Career Grant CUHK24300915. 
This collaboration was initiated at the Hausdorff Center for Mathematics and Oberwolfach, and continued at the joint AMS-EMS-SPM 2015 international meeting at Porto. The authors thank all institutions involved for gracious and productive work environments.

\bibliographystyle{alpha}
\bibliography{GPRY}
\vspace{.25cm}
{\tiny
\noindent \textsc{Shaoming Guo, Indiana University Bloomington, 107 S Indiana Ave, Bloomington, IN 47405, USA}\\
{\em email: }\textsf{\bf shaoguo@iu.edu}\\

\noindent \textsc{Lillian B. Pierce, Duke University, Box 90320, 120 Science Drive, Durham NC 27708, USA}\\
{\em email: }\textsf{\bf pierce@math.duke.edu}\\

\noindent \textsc{Joris Roos, University of Bonn, Mathematical Institute, Endenicher Allee 60, 53115 Bonn, Germany}\\
{\em email: }\textsf{\bf jroos@math.uni-bonn.de}\\

\noindent \textsc{Po-Lam Yung, The Chinese University of Hong Kong, Ma Liu Shui, Shatin, Hong Kong}\\
{\em email: }\textsf{\bf plyung@math.cuhk.edu.hk}
}
\end{document}